\documentclass[reqno]{amsart}
\numberwithin{equation}{section}
\usepackage[english]{babel}
\usepackage[T1]{fontenc}
\usepackage[latin1]{inputenc}
\usepackage{indentfirst}
\usepackage{amsmath,amssymb, amsbsy}
\usepackage{amsfonts}
\usepackage{color}
\usepackage{latexsym}
\usepackage{eucal}
\usepackage{amsthm}
\usepackage[dvips]{graphicx}
\DeclareGraphicsExtensions{.pdf,.png,.jpg,.eps}

\usepackage{color}
\usepackage[latin1]{inputenc}
\usepackage[active]{srcltx}
\newtheorem{teo}{Theorem}[section]

\newtheorem{Lemma}[teo]{Lemma}
\newtheorem{Proposition}[teo]{Proposition}
\theoremstyle{definition}
\newtheorem{Definition}[teo]{Definition}

\begin{document}

\author{Susanna Terracini and Stefano Vita}

\address{Susanna Terracini and Stefano Vita \newline \indent
 Dipartimento di Matematica ``Giuseppe Peano'', \newline \indent
Universit\`a di Torino, \newline \indent
Via Carlo Alberto, 10,
10123 Torino, Italy}
\email{susanna.terracini@unito.it}
\email{{\tt stefano.vita@unito.it}}


\title[Growth of competing densities for a nonlocal elliptic blow-up system]
{On the asymptotic growth of positive solutions to a nonlocal elliptic blow-up system involving strong competition}
\date{\today}


\thanks{Work partially supported by the 
ERC Advanced Grant 2013 n. 339958
{\it Complex Patterns for Strongly Interacting Dynamical Systems - COMPAT} and by the PRIN-2012-74FYK7 Grant {\it Variational and perturbative aspects of nonlinear differential problems}. The authors wish to thank Alessandro Zilio for many fruitful conversations.
}

\maketitle


\begin{abstract}
For a competition-diffusion system involving the fractional Laplacian of the form
\begin{equation*}\label{syst1}
-(-\Delta)^su=uv^2,\quad-(-\Delta)^sv=vu^2,\quad u,v>0 \ \mathrm{in} \ \mathbb{R}^N,
\end{equation*}
whith $s\in(0,1)$, we prove that the maximal asymptotic growth rate for its entire solutions is $2s$. Moreover, since we are able to construct symmetric solutions to the problem, when $N=2$ with prescribed growth arbitrarily close to the critical one, we can conclude that the asymptotic bound found is optimal.  Finally, we prove existence of genuinely higher dimensional solutions, when $N\geq 3$. Such problems arise, for example, as blow-ups of fractional reaction-diffusion  systems when the interspecific competition rate tends to infinity. 
\end{abstract}

\section{Introduction and main results}
This paper deals with the existence and classification of positive entire solutions to polynomial systems involving the (possibly) $s$-fractional Laplacian of the following form:
\begin{equation*}\label{syst1}
-(-\Delta)^su=uv^2,\quad-(-\Delta)^sv=vu^2,\quad u,v>0 \ \mathrm{in} \ \mathbb{R}^N.
\end{equation*}
Such systems arise,  for example, as blow-ups of fractional reaction-diffusion systems when the interspecific competition rate tends to infinity. In this framework, the existence and classification of entire solutions plays a key role in the asymptotic analysis (see, for instance, \cite{soazil2015, soazil2017}). The case of standard diffusion ($s=1$) has been intensively treated in the recent literature, also in connection with a De Giorgi-like conjecture about monotone solutions being one dimensional. In particular, a complete classification of solutions having linear (the lowest possible growth rate) has been given in \cite{berelin, bereterra, farina,farsoa,soazil1,wang}. On the other hand, when $s=1$, positive solutions having arbitrarily large polynomial growth were discovered in \cite{bereterra} and with exponential growth in \cite{soazil2014}.\medskip

Competition-diffusion nonlinear systems with $k$-components involving the fractional Laplacian have been the object of a recent literature, starting with  \cite{terverzil,terverzil2}, where the authors
provided asymptotic estimates for solutions to systems 
of the form
\begin{equation}\label{systgoss}
\begin{cases}
(-\Delta)^su_i=f_{i,\beta}(u_i)-\beta u_i\sum_{j\neq i}a_{ij}u_j^2, &i=1,...,k,\\
u_i\in H^s(\mathbb{R}^N),
\end{cases}
\end{equation}
where $N\geq 2$, $a_{ij}=a_{ji}>0$, when $\beta>0$ (the competition parameter) goes to $+\infty$. Moreover we consider $f_{i,\beta}$ as continuous functions which are uniformly bounded on bounded sets with respect to $\beta$ (see \cite{terverzil,terverzil2} for details). The fractional Laplacian is defined for every $s\in(0,1)$ as
\begin{equation*}\label{fraclapl}
(-\Delta)^su(x)=c(N,s) \ \mathrm{PV}\int_{\mathbb{R}^N}\frac{u(x)-u(y)}{|x-y|^{N+2s}} \ \mathrm{d}y.
\end{equation*}
In order to state our results, we adopt the approach  of  Caffarelli-Silvestre \cite{cafsil}, and we see the fractional Laplacian as a Dirichlet-to-Neumann operator; that is, we consider the  extension problem for \eqref{systgoss}. In other words, we study an auxiliary problem in the upper half space in one more dimension\footnote{Throughout this paper we assume the following notations: $z=(x,y)$ denotes a point in $\mathbb{R}^{N+1}_+$, with $x\in\partial\mathbb{R}^{N+1}_+:=\mathbb{R}^N$ and $y\in\mathbb{R}_+$. Moreover, $B^{+}_r(z_0):=B_r(z_0)\cap\mathbb{R}^{N+1}_+$ is the half ball, and its boundary is divided in the hemisphere $\partial^+ B^{+}_r(z_0):=\partial B^{+}_r(z_0)\cap\mathbb{R}^{N+1}_+$ and in the flat part $\partial^0 B^{+}_r(z_0):=\partial B^{+}_r(z_0)\setminus\partial^+ B^{+}_r(z_0)$. When the center of balls and spheres is omitted 
, then $z_0=0$.}; that is, letting $a:=1-2s$, for any $i=1,...,k$ the localized version of \eqref{systgoss},
\begin{equation}\label{cafsilvess1}
\begin{cases}
L_au_i=0, &\mathrm{in} \ B^+_1\subset\mathbb{R}^{N+1}_+,\\
\partial^a_yu_i=f_{i,\beta}(u_i)-\beta u_i\sum_{j\neq i}a_{ij}u_j^2,  &\mathrm{in} \ \partial^0B^+_1\subset\partial\mathbb{R}^{N+1}_+=\mathbb{R}^{N}\times\{0\},
\end{cases}
\end{equation}
where the degenerate/singular elliptic operator $L_a$ is defined as
\begin{equation*}\label{Lau}
-L_au:=\mathrm{div}(y^a\nabla u),
\end{equation*}
and the linear operator $\partial^a_y$ is defined as
\begin{equation*}\label{pary}
-\partial^a_yu:=\lim_{y\to 0^+}y^a\frac{\partial u}{\partial y}.
\end{equation*}
The new problem \eqref{cafsilvess1}  is equivalent to the original when we deal with solutions in the energy space associated with the two operators. In fact a solution $U$ to the extension problem is the extension of the correspondent solution $u$ of the original nonlocal problem in the sense that $U(x,0)=u(x)$. Let us remark that if $s=\frac{1}{2}$, then $a=0$ and hence $L_0=-\Delta$ and the boundary operator $-\partial^0_y$ becomes the usual normal derivative $\partial_y$. Moreover we remark that the extension problem has a variational nature in some weighted Sobolev spaces related to the Muckenhoupt $A_2$-weights (see for instance \cite{kufner}). Hence, given $\Omega\subset\mathbb{R}^{N+1}_+$, we can introduce the Hilbert spaces
\begin{equation*}
H^{1;a}(\Omega):=\left\{ u:\Omega\to\mathbb{R} \ : \ \int_{\Omega}y^a(|u|^2+|\nabla u|^2)<+\infty \right\},
\end{equation*}
and
\begin{equation*}
H^{1;a}_{loc}\left(\overline{\mathbb{R}^{N+1}_+}\right):=\left\{ u:\mathbb{R}^{N+1}_+\to\mathbb{R} \ : \ \forall r>0, u|_{B_r^+}\in H^{1;a}(B_r^+) \right\},
\end{equation*}
where the functions $u=u(z)$ are functions of the variables $z=(x,y)\in\mathbb{R}^{N+1}_+$. In the quoted papers \cite{terverzil,terverzil2},  the authors make use of Almgren's and Alt-Caffarelli-Friedman's type monotonicity formul\ae \ in order to obtain uniform H\"older bounds with small exponent $\alpha=\alpha(N,s)$ for bounded energy solutions of the Gross-Pitaevskii system. Passing to the limit as the competition parameter $\beta\longrightarrow+\infty$ and using suitably rescaled dependent and independent variables in \eqref{cafsilvess1},  a main step consists in classifying the entire solutions to the limiting system solved by blow-up solutions. In particular, we are interested in studying some qualitative properties related to the asymptotic growth for positive entire solutions of this elliptic system in case of two components. In our setting, the resulting system is the following
\begin{equation}\label{cafsilvess}
\begin{cases}
L_au=L_av=0, &\mathrm{in} \ \mathbb{R}^{N+1}_+,\\
u,v>0, &\mathrm{in} \ \overline{\mathbb{R}^{N+1}_+},\\
-\partial^a_yu=uv^2, \ -\partial^a_yv=vu^2, &\mathrm{in} \ \partial\mathbb{R}^{N+1}_+,
\end{cases}
\end{equation}
which is equivalent to
\begin{equation}\label{syst0}
-(-\Delta)^{s}u=uv^2,\quad-(-\Delta)^{s}v=vu^2,\quad u,v>0 \ \mathrm{in} \ \mathbb{R}^N.
\end{equation}
We focus our attention on positive solutions since this condition follows requiring that the original Gross-Pitaevskii solutions do not change sign in $\mathbb{R}^N$. Some relevant qualitative properties of positive solutions to system \eqref{syst0} have been recently investigated by Wang and Wei in \cite{wangwei}. In particular, they proved uniqueness for the one-dimensional solutions when $s>1/4$, up to translation and scaling. Moreover, they highlighted a universal polynomial bound at infinity for positive subsolutions. Their result shows a striking contrast between the cases of the fractional and the local diffusion; indeed, in the latter case, there are solutions  having arbitrarily large polynomial and even exponential  growth \cite{bereterra,soazil2014}. As the polynomial bound in \cite{wangwei} is restricted to positive solutions and there are sign-changing solutions to the equation $L_a u=0$ having arbitrarily large growth rate, we suggest that the picture may change also considering sign-changing solutions to the Gross-Pitaevskii system. \medskip

Following \cite{soater}, we give the following definition.
\begin{Definition}
Let $(u,v)$ be a solution to \eqref{cafsilvess}. We say that $(u,v)$ has algebraic growth if there exist two constants $c,d>0$ such that
\begin{equation}\label{growthalgeb}
u(x,y)+v(x,y)\leq c\left(1+|x|^2+y^2\right)^{d/2}\qquad\forall (x,y)\in\overline{\mathbb{R}^{N+1}_+}.
\end{equation}
Moreover we say that $(u,v)$ has growth rate $d>0$ if
\begin{equation}\label{eq:whosd}
\lim_{r\to+\infty}\dfrac{\int_{\partial^+B_r^+}y^a(u^2+v^2)}{r^{N+a+2d'}}=\begin{cases}
+\infty & \mathrm{if} \ d'<d\\
0 & \mathrm{if} \ d'>d.
\end{cases}
\end{equation}
\end{Definition}
It can be shown that the threshold exponent $d$ appearing in \eqref{eq:whosd} is exactly the extremal one for which \eqref{growthalgeb} holds (see Proposition \ref{prop:growthrate}).\medskip

The aim of our work is to find the maximal asymptotic growth for positive solutions to \eqref{syst0}; to this aim, we shall construct a family of solutions possessing some natural symmetry, this extending the results of \cite{bereterra} to the case of fractional diffusions.
\medskip

In what follows, we will study an eigenvalue problem for the spherical part of the operator $L_a$. We can think to such a operator as a Laplace-Beltrami-type operator on the superior hemisphere $S^N_+$ of the unit sphere $S^N\subset\mathbb{R}^{N+1}$. Our aim is to deal with some $\mathbb{G}_k$-equivariant optimal partitions, in the case $N=2$, where the symmetry group $\mathbb{G}_k$ acts cyclically with order $k$. In particular, we will construct a sequence of optimal partition first-eigenvalues $\{\lambda_1^s(k)\}_{k=1}^{+\infty}$ and related nonnegative eigenfunctions $\{u_k\}_{k=1}^{+\infty}$, where $k$ is the order of the symmetry group imposed on the boundary condition region.\medskip

Hence we will prove the following asymptotic bound.
\begin{teo}\label{tteo1}
Let $s\in(0,1)$ and $N\geq 2$. Let $(u,v)$ be a positive solution to \eqref{cafsilvess}. Then, there exists a constant $c>0$ such that
\begin{equation}\label{lineargrowth}
u(x,y)+v(x,y)\leq c\left(1+|x|^2+y^2\right)^{s}.
\end{equation}
\end{teo}
Hence, we will use the sequence of eigenfunctions previously seen, in order to construct a sequence of positive solutions to \eqref{cafsilvess} possessing some symmetries and having an asymptotic growth rate arbitrarily close to to the critical one; that is, we will prove
\begin{teo}\label{tteo2}
When $N=2$ and $s\in(0,1)$ there exists a sequence of positive solutions $(u_k,v_k)$ to the system \eqref{cafsilvess} having growth rate $d(k)\in [s,2s)$, where $d(k)$ converges monotonically to $2s$.
\end{teo}
These prescribed growth solutions for \eqref{cafsilvess} in space dimension $N=2$ are also solutions with the same properties for the same problem in any higher dimension.\medskip

Eventually, in the last section, we will show the existence of entire solutions to \eqref{cafsilvess} which are truly $N$-dimensional, in the sense that they can not be obtained by adding coordinates in a constant way starting from a $2$-dimensional solution.

\section{Bound on the growth rate of positive solutions}
Our first general purpose is to study the asymptotic behavior of entire nonnegative solutions to the cubic system
\begin{equation*}\label{systxx}
-(-\Delta)^{s}u=uv^2,\quad-(-\Delta)^{s}v=vu^2,\quad u,v>0 \ \mathrm{in} \ \mathbb{R}^N.
\end{equation*}
In particular we prove that solutions can not grow faster than $2s$ at infinity. Furthermore, as we will are able to construct solutions to this problem with prescribed growth rate arbitrarily close to the critical one, we can conclude that this asymptotic bound is optimal. As said in the introduction, we will deal with the equivalent Caffarelli-Silvestre extension problem defined in \eqref{cafsilvess}.

First we will introduce the Almgren frequency function and its monotonicity formula which are the main instruments that we need to prove Theorem \ref{tteo1} and Theorem \ref{tteo2}.

\subsection{Almgren monotonicity formula}
Now, we are going to summarize some results proved in \cite{terverzil,terverzil2,wangwei}, involving the Almgren monotonicity formula for solutions to \eqref{cafsilvess}. First, solutions of \eqref{cafsilvess} satisfy a Pohozaev identity; that is, for any $x_0\in\mathbb{R}^N$ and $r>0$,
\begin{eqnarray}
(N-1+a)\int_{B_r^+(x_0,0)}y^a(|\nabla u|^2+|\nabla v|^2)&=&r\int_{\partial^+B_r^+(x_0,0)}y^a(|\nabla u|^2+|\nabla v|^2)-2y^a(|\partial_ru|^2+|\partial_rv|^2)\nonumber\\
&+&r\int_{S^{N-1}_r(x_0,0)}u^2v^2-N\int_{\partial^0B_r^+(x_0,0)}u^2v^2.
\end{eqnarray}
Moreover, let us recall the following definitions
\begin{equation}\label{E}
E(r,x_0;u,v):=\frac{1}{r^{N-1+a}}\left(\int_{B_r^+(x_0,0)}y^a(|\nabla u|^2+|\nabla v|^2) \ + \ \int_{\partial^0B_r^+(x_0,0)}u^2v^2\right),
\end{equation}
and
\begin{equation}\label{H}
H(r,x_0;u,v):=\frac{1}{r^{N+a}}\int_{\partial^+B_r^+(x_0,0)}y^a(u^2+v^2).
\end{equation}
Hence, defining the frequency as $N(r,x_0;u,v):=\frac{E(r,x_0;u,v)}{H(r,x_0;u,v)}$, the Almgren monotonicity formula holds; that is, the frequency $N(r,x_0;u,v)$ is non decreasing in $r>0$. Moreover, if $(u,v)$ is a solution to \eqref{cafsilvess} and $N(R)\geq d$ then for $r>R$ it holds that $H(r)/r^{2d}$ is non decreasing in $r$. Hence, if we consider $(u,v)$ a solution of \eqref{cafsilvess} on a bounded half ball $B^+_R$ and if $N(R)\leq d$, then for every $0<r_1\leq r_2\leq R$ it holds that
\begin{equation}\label{rappHH}
\frac{H(r_2)}{H(r_1)}\leq e^{\frac{d}{1-a}}\frac{r_2^{2d}}{r_1^{2d}}.
\end{equation}

\subsection{Eigenvalue problem for a Laplace-Beltrami-type operator with mixed boundary conditions}
As the authors of \cite{terverzil,terverzil2,wangwei} have pointed out, the regularity and the asymptotic growth of solutions to competition problems are related to an optimal partition problem on the superior hemisphere $S^N_+\subset\mathbb{R}^{N+1}_+$. Likewise the case of the Laplacian, we wish to express the extension operator $L_a$ in spherical coordinates, in order to write it as the sum of a radial part and a Laplace-Beltrami-type operator defined on the superior hemisphere (see \cite{ruland}). Let us consider in $\mathbb{R}^{N+1}_+$ the spherical coordinates $(r,\theta,\phi)$ such that $y=r\sin\theta$, with $\theta\in[0,\pi/2]$ and $\phi=(\phi_1,...,\phi_{N-1})$ parametrizing the position over $S^{N-1}\subset\mathbb{R}^N$. Hence,
\begin{equation}\label{decomp}
-L_au=\nabla\cdot y^a\nabla u=(\sin\theta)^a\frac{1}{r^N}\partial_r(r^{N+a}\partial_ru)+r^{a-2}L_a^{S^N}u,
\end{equation}
where the Laplace-Beltrami-type operator is defined as
\begin{equation}
L_a^{S^N}u:=\nabla_{S^N}\cdot(\sin\theta)^a\nabla_{S^N}u=\nabla_{S^N}\cdot y^a\nabla_{S^N}u,
\end{equation}
and $\nabla_{S^N}$ is the tangential gradient on $S^N_+$. For every open $\omega\subset S^{N-1}:=\partial S^N_+$, we define the first $s$-eigenvalue associated to $\omega$ as 
\begin{equation}\label{lambda1somega}
\lambda_1^s(\omega):=\inf\left\{\frac{\int_{S^N_+}y^a|\nabla_{S^N}u|^2}{\int_{S^N_+}y^au^2} : u\in H^{1;a}(S^N_+)\setminus\{0\}, \ u=0 \mathrm{ \ in \ } S^{N-1}\setminus\omega\right\}.
\end{equation}
So, such a minimization problem has a natural variational structure on the weighted Sobolev space $H^{1;a}(S^N_+):=\left\{u:S^N_+\to\mathbb{R} \ : \ \int_{S^N_+}y^a|\nabla_{S^N}u|^2+\int_{S^N_+}y^au^2<+\infty\right\}$; which is an Hilbert space. In fact, defining $H^{1;a}_{\omega}(S^N_+):=\{u\in H^{1;a}(S^N_+) : \ u=0 \mathrm{ \ in \ } S^{N-1}\setminus\omega\}$ for every fixed $\omega\subset S^{N-1}$, we get in this space the existence of a nontrivial and nonnegative minimizer of the Rayleigh quotient
\begin{equation*}
\mathcal{R}^a(u):=\frac{\int_{S^N_+}y^a|\nabla_{S^N}u|^2}{\int_{S^N_+}y^au^2},
\end{equation*}
which is also an eigenfunction related to $\lambda^s_1(\omega)$ since it is a weak solution to the following mixed Dirichlet-to-Neumann boundary eigenvalue problem for the spherical part of the $L_a$ operator
\begin{equation}\label{Plambdagamma}
\begin{cases}
-L_a^{S^N} u=y^a\lambda_1^s(\omega) u  &\mathrm{in}\quad S^N_+, \\
u=0  &\mathrm{in}\quad S^{N-1}\setminus\omega, \\
\partial^a_y u=0 &\mathrm{in}\quad\omega\subset S^{N-1}.
\end{cases}
\end{equation}
Moreover, for every $\omega\subset S^{N-1}$ it holds that
\begin{equation*}
H^{1;a}_{0}(S^N_+)\subseteq H^{1;a}_{\omega}(S^N_+)\subseteq H^{1;a}(S^N_+).
\end{equation*}
Hence by definition, for any  $\omega\subset S^{N-1}$,
\begin{equation}\label{zeig}
\lambda_1^s(S^{N-1})\leq\lambda_1^s(\omega)\leq\lambda_1^s(\emptyset).
\end{equation}
Let us now define the characteristic exponent
\begin{equation}\label{mono}
d(t):=\sqrt{\left(\frac{N-2s}{2}\right)^2+t}-\frac{N-2s}{2}.
\end{equation}
The characteristic exponent is defined in such a way that $u$ is a nonnegative eigenfunction of $\lambda_1^s(\omega)$ if and only if its $d(\lambda_1^s(\omega))$-homogeneous extension to $\mathbb{R}^{N+1}_+$ is $L_a$-harmonic.

Let us define by $\omega^c=S^{N-1}\setminus\overline{\omega}$, with $\omega\subset S^{N-1}$ open. Obviously $\omega\cap\omega^c=\emptyset$ and $\overline{\omega}\cup\overline{\omega^c}=S^{N-1}$. From now on, we suppose that $\gamma=\overline{\omega}\cap\overline{\omega^c}$ is a $(N-2)$-dimensional smooth submanifold. Analogously with the case of the Laplacian  in \cite{coloperal}, one can consider two nonnegative eigenfunctions $u_1,u_2$ of \eqref{Plambdagamma} with eigenvalues $\lambda_1^s(\omega_1)$ and $\lambda_1^s(\omega_2)$. In our setting, if there exists $\alpha\in(0,1)$ such that $u_1,u_2\in C^{0;\alpha}(\overline{S^N_+})$, $\mathcal{H}^{N-1}(\omega_1)>\mathcal{H}^{N-1}(\omega_2)$ and $\omega_2\subset\omega_1$, then it holds that
\begin{equation}\label{coloper}
\lambda_1^s(\omega_1)<\lambda_1^s(\omega_2).
\end{equation}
In fact, integrating by parts with respect to both the eigenfunctions the quantity $\int_{S^N_+}y^a\nabla_{S^N}u_1\nabla_{S^N}u_2$, we find
\begin{equation}
\lambda_1^s(\omega_1)\int_{S^N_+}y^au_1u_2+\int_{\omega_1^c\cap\omega_2}(\partial_y^au_1)u_2=\lambda_1^s(\omega_2)\int_{S^N_+}y^au_1u_2+\int_{\omega_2^c\cap\omega_1}(\partial_y^au_2)u_1,
\end{equation}
and since $\omega_2\subset\omega_1$, then $\omega_1^c\cap\omega_2=\emptyset$ and $\omega_2^c\cap\omega_1=\omega_3$ open. Hence, \eqref{coloper} holds using the Hopf lemma
\begin{equation}
(\lambda_1^s(\omega_1)-\lambda_1^s(\omega_2))\int_{S^N_+}y^au_1u_2=\int_{\omega_3}(\partial_y^au_2)u_1<0.
\end{equation}

\subsection{Blow-down analysis and the maximal growth rate}
Now, after performing a scaling (blow-down) analysis over general positive solutions to \eqref{cafsilvess}, we will prove the upper bound on the growth at infinity; that is, Theorem \ref{tteo1}. First, we summarize the steps done by Wang and Wei. Theorem 2.3 in \cite{wangwei} proves that, taking a positive solution $(u,v)$ to \eqref{cafsilvess}, then there exist two constants $d,c>0$ such that
\begin{equation}\label{dgrowth}
u(x,y)+v(x,y)\leq c\left(1+|x|^2+y^2\right)^{d/2}.
\end{equation}
Moreover, in Proposition 3.5, they proved that condition \eqref{dgrowth} is equivalent to the following upper bound over the frequency
\begin{equation}\label{prop3.5}
N(R)\leq d, \quad\forall R>0.
\end{equation}
We can consider $d>0$ which is the infimum such that condition \eqref{dgrowth} holds. For such a number, if there exists the limit $\lim_{R\to+\infty}N(R)$, then of course it is exactly equal to $d$. In other words, we have:
\begin{Proposition}\label{prop:growthrate}
The growth rate of a positive solution $(u,v)$ to \eqref{cafsilvess} is $d$ if and only if 
\[\lim_{R\to+\infty}N(R)=d\;,
\]

\end{Proposition}

\subsubsection{Proof of Theorem \ref{tteo1}}
Let $(u,v)$ be a positive solution to \eqref{cafsilvess}. Note that \eqref{prop3.5} combined with the Almgren monotonicity formula also implies that $\lim_{R\to+\infty}N(R)=d$. Let us define for $R\longrightarrow+\infty$ the blow-down sequence
\begin{equation*}\label{bds}
u_R(z):=L(R)^{-1}u(Rz),\quad v_R(z):=L(R)^{-1}v(Rz),
\end{equation*}
with $L(R)$ taken so that $H((u_R,v_R),1)=1$. So, the sequence satisfies
\begin{equation*}\label{sblowk}
\begin{cases}
L_a u_R=L_a v_R=0  & \mathrm{in}\quad\mathbb{R}^{N+1}_+,\\
-\partial^a_y u_R=\kappa_Ru_Rv_R^2, \ -\partial^a_y v_R=\kappa_Rv_Ru_R^2 & \mathrm{in}\quad\partial\mathbb{R}^{N+1}_+,
\end{cases}
\end{equation*}
where $\kappa_R=L(R)^2R^{1-a}$. By the Liouville theorem (see Proposition 3.9 in \cite{terverzil2}), for some $\alpha>0$ small there exists a constant $C_\alpha$ such that $L(R)\geq C_\alpha R^\alpha$ so that $\kappa_R\longrightarrow+\infty$ as $R\longrightarrow+\infty$. Hence, thanks to \eqref{rappHH} we get the following integral uniform upper bound; that is, $H((u_R,v_R),r)\leq r^{2d}$ for every $r>1$. Since $(u_R,v_R)$ satisfy the requirements of Lemma A.2 in \cite{wangwei}, for every $r>1$ we get that
\begin{equation*}\label{uppp}
\sup_{B^+_r}(u_R+v_R)\leq Cr^d.
\end{equation*}
Then, thanks to the uniform H\"older estimates proved in \cite{terverzil2}, for some small $\alpha>0$, the sequence $\{(u_R,v_R)\}$ is uniformly bounded in $C_{loc}^{0,\alpha}(\overline{\mathbb{R}^{N+1}_+})$. Hence, letting $R\longrightarrow+\infty$, up to consider a subsequence, we get weakly convergence in $H^{1;a}_{loc}(\mathbb{R}^{N+1}_+)$ and uniform convergence in $C_{loc}^{0,\alpha}(\overline{\mathbb{R}^{N+1}_+})$ of the sequence $\{(u_R,v_R)\}$ to a couple of functions $(u_{\infty},v_{\infty})$ which are segregated in $\partial\mathbb{R}^{N+1}$ in the sense that $u_{\infty}v_{\infty}=0$ in $\partial\mathbb{R}^{N+1}$. Proceeding as in \cite{wangwei}, using the fact that $N((u_{\infty},v_{\infty}),r)=d$ for any $r>0$, we can conclude that such functions are homogeneous of degree $d$ and segregated in $\partial\mathbb{R}^{N+1}_+$; that is, they solve the following problem
\begin{equation}\label{sblowcseg}
\begin{cases}
L_a u_{\infty}=L_a v_{\infty}=0 & \mathrm{in}\quad\mathbb{R}^{N+1}_+, \\
u_{\infty}\partial^a_yu_{\infty}=v_{\infty}\partial^a_yv_{\infty}=0 & \mathrm{in}\quad\partial\mathbb{R}^{N+1}_+, \\
u_{\infty}v_{\infty}=0 & \mathrm{in}\quad\partial\mathbb{R}^{N+1}_+.
\end{cases}
\end{equation}
Moreover, such solutions have the form
\begin{equation*}\label{uinf}
u_{\infty}(r,\theta)=r^dg(\theta),\quad v_{\infty}(r,\theta)=r^dh(\theta),
\end{equation*}
where $g,h$ are defined on the upper hemisphere $S^N_+=\partial^+B^+_1$. Since we have constructed the blow-down sequence so that $H((u_R,v_R),1)=1$, then
\begin{equation}
\int_{S^N_+}y^a(g^2+h^2)=1,
\end{equation}
and hence can not happen that both $g$ and $h$ vanish identically in $S^N_+$, but at most only one component is identically zero. In any case, by the homogeneity of the blow-down limit and the fact that $(u_{\infty},v_{\infty})$ are $L_a$-harmonic, any nontrivial component is an eigenfunction for the spherical part of $L_a$ in the sense seen in \eqref{Plambdagamma} on $S^N_+$. Moreover, such eigenfunction must own eigenvalue $\lambda$ which has the following relation with the characteristic exponent $d$,
\begin{equation}\label{lambdad}
\lambda=d(d+N-1+a).
\end{equation}
But we have seen with \eqref{zeig} that such eigenvalue can not be larger than $\lambda_1^s(\emptyset)$, achieved by $u(x,y)=y^{2s}$ which has $d(\lambda_1^s(\emptyset))=2s$. Moreover, by \eqref{mono}, the map $t\longmapsto d(t)$ is strictly increasing and hence $d\leq d(\lambda_1^s(\emptyset))$. By \eqref{dgrowth}, Theorem \ref{tteo1} is proved.

\section{Prescribed growth solutions}
From now on in this section we consider the case $N=2$ and we study the optimal boundary condition minimizing the first eigenvalue of \eqref{Plambdagamma} under some requirements over the measure and the symmetries of $\omega\subset S^1$. Doing this, we will be able to construct positive solutions to \eqref{cafsilvess} with prescribed growth and depending in some way on the 2-dimensional eigenvalue problem.

In the next section, we are going to introduce a suitable type of Schwarz symmetrization, that will be the main tool that we need to study this optimal boundary condition problem.

\subsection{Polarization and foliated Schwarz symmetrization}
From now on we follow some ideas contained in \cite{brocksolynin,smetswillem}. We can state the results in this section in any space dimension $N\geq 2$. Let us define by $\mathcal{H}$ the set of all half spaces in $\mathbb{R}^{N+1}$ determined by the set of all the affine hyperplanes with orientation, and by $\mathcal{H}_0$ the subset of $\mathcal{H}$ determined by the euclidean hyperplanes with orientation. Let $H\in\mathcal{H}$ be a half space, we denote by $\sigma_H$ the reflection with respect to the hyperplane $\partial H$.
\begin{Definition}\label{def1}
Let $H\in\mathcal{H}$ be a half space. The polarization of a measurable nonnegative function $u$ with respect to $H$ is the function defined by
\begin{equation*}\label{pola}
u_H(z):=\begin{cases}
\max\{u(z),u(\sigma_H(z))\} &\mathrm{if} \ z\in H,\\
\min\{u(z),u(\sigma_H(z))\} &\mathrm{if} \ z\in \mathbb{R}^{N+1}\setminus H.
\end{cases}
\end{equation*}
\end{Definition}
In the same way we can define the polarization $A_H$ of a set $A\subset\mathbb{R}^{N+1}$ with respect to $H\in\mathcal{H}$ in the sense that $\chi_{A_H}=(\chi_{A})_H$. It is well known that the polarization mapping $A\mapsto A_H$ is a rearrangement of $\mathbb{R}^{N+1}$ for the Lebesgue measure for any $H\in\mathcal{H}$; that is, it satisfies both the monotonicity property ($A\subset B\Rightarrow A_H\subset B_H$) and the measure conservation property ($\mathcal{L}^{N+1}(A_H)=\mathcal{L}^{N+1}(A)$) (see \cite{smetswillem}).

Let us consider $\Sigma_1=\{x_1=0\}$ as a fixed hyperplane ($\Sigma_1=\partial H_1$ with $H_1=\{x_1>0\}$), and denote by $\sigma_1:=\sigma_{\Sigma_1}$ the reflection with respect to $\Sigma_1$. Let us now consider the point $z_0^1\in\overline{S^N_+}$ which maximizes the distance from the hyperplane $\Sigma_1$ (actually, there are two points with this property $z_0^1,z_0^2$, we choose the one in $H_1$). This point lies on $S^{N-1}=\partial S^N_+$. Let us define $\mathcal{H}_1:=\{H\in\mathcal{H}_0: \ z_0^1\in H \mathrm{ \ and \ axis \ } y \mathrm{ \ lies \ on \ }\partial H\} $. Since the measure given by $\mathrm{d}\mu:=y^a\mathrm{d}S_N(z)$ is mapped into itself by the reflection $\sigma_H$ for any $H\in\mathcal{H}_1$, with the same arguments in \cite{smetswillem}, we can see that polarization is also a rearrangement of $S^N_+$ for the measure $\mu$ for any $H\in\mathcal{H}_1$. Moreover, we can obtain the invariance of the norm in weighted spaces under polarization for $H\in\mathcal{H}_1$; that is, when $u\in L^p(S^N_+;\mathrm{d}\mu)$ with $1\leq p<+\infty$, we have $u_H\in L^p(S^N_+;\mathrm{d}\mu)$ with
\begin{equation}\label{normwh}
\int_{S^N_+}y^a|u_H|^p\mathrm{d}S_N=\int_{S^N_+}y^a|u|^p\mathrm{d}S_N,
\end{equation}
and if $u\in W^{1,p}_+(S^N_+;\mathrm{d}\mu)$ with $1\leq p<+\infty$, hence $u_H\in W^{1,p}_+(S^N_+;\mathrm{d}\mu)$ with
\begin{equation}\label{normlh}
\int_{S^N_+}y^a|\nabla_{S^N}u_H|^p\mathrm{d}S_N=\int_{S^N_+}y^a|\nabla_{S^N}u|^p\mathrm{d}S_N.
\end{equation}
Now we want to define the foliated Schwarz symmetrization on the hemisphere. Consider for $\overline{y}\in[0,1)$ the $(N-1)$-sphere defined by
\begin{equation*}\label{steta}
S^{N-1}_{\overline{y}}:=\overline{S^N_+}\cap\{y=\overline{y}\}.
\end{equation*}
Let us define on every $(N-1)$-sphere $S^{N-1}_{y}$ the point $z^1_{y}$ so that it has the same parametrizing angle $\overline\phi$ of the point $z^1_0$. The symmetrization $A^*$ of a set $A\subset S^{N-1}_{y}$ with respect to $z^1_{y}$ is defined as the closed geodesic ball centered in $z^1_{y}$ such that $\mathcal{L}^{N-1}(A^*)=\mathcal{L}^{N-1}(A)$. The symmetric decreasing rearrangement $f^*$ of a nonnegative measurable function $f$ defined on $S^{N-1}_{y}$ is such that $\{f>t\}^*=\{f^*>t\}$ for every $t\geq 0$. We remark that this symmetrization is a rearrangement of the sphere $S^{N-1}_y$ for the measure $\mathcal{L}^{N-1}$, for every fixed $y\in[0,1)$.
\begin{Definition}\label{def2}
Let $u\in H^{1;a}(S^N_+)$ be a nonnegative function. The foliated Schwarz symmetrization $u^*$ of $u$ is defined on the hemisphere $S^N_+$ by the symmetric decreasing rearrangement of the restriction of $u$ on every $S^{N-1}_{y}$; that is, $u^*|_{S^{N-1}_{y}}=(u|_{S^{N-1}_{y}})^*$ for every $y\in[0,1)$.
\end{Definition}
One can check that also the foliated Schwarz symmetrization is a rearrangement of $S^N_+$ for $\mu$, since it satisfies both the monotonicity property ($A\subset B\Rightarrow A^*\subset B^*$) and the measure conservation property ($\mu(A^*)=\mu(A)$), where the symmetrization $A^*$ of a set $A\subset S^N_+$ is defined as the only set in $S^N_+$ such that $A^*\cap S^{N-1}_{y}=(A\cap S^{N-1}_{y})^*$ for every $y\in[0,1)$, in the sense of symmetrization of a set in $S^{N-1}_{y}$ given previously (the idea is that this symmetrization map works only on the $x$-variable and so $\mathrm{d}\mu$ is mapped into itself). Moreover, it is easy to see that for every nonnegative $u\in H^{1;a}(S^N_+)$ and for every $H\in\mathcal{H}_1$ it holds that
\begin{equation}\label{invariance}
(u^*)_H=u^*=(u_H)^*.
\end{equation}
Hence it holds the following result from \cite{smetswillem}. For completeness we adapt to our hemispherical case the proof of Smets and Willem.
\begin{Lemma}\label{lem1}
Let $u\in C(S^N_+)$ be a nonnegative function. If $u\neq u^*$, then there exists $H\in\mathcal{H}_1$ such that
\begin{equation}\label{strictt}
||u_H-u^*||_{L^{2}(S^N_+;\mathrm{d}\mu)}<||u-u^*||_{L^{2}(S^N_+;\mathrm{d}\mu)}.
\end{equation}
\end{Lemma}
\proof
First of all, we remark that always the non strict inequality in \eqref{strictt} holds (rearrangement for a suitable measure $\mu$ is a contraction in $L^p(\mathrm{d}\mu)$ for any $1\leq p<+\infty$). If $u\neq  u^*$, there exists $y\in[0,1)$ and $t\geq 0$ such that $\{u>t\}\cap S^{N-1}_{y}\neq \{u^*>t\}\cap S^{N-1}_{y}$ and since the foliated Schwarz symmetrization is a rearrangement, then $\mathcal{L}^{N-1}(\{u>t\}\cap S^{N-1}_{y})=\mathcal{L}^{N-1}(\{u^*>t\}\cap S^{N-1}_{y})$; so, by the continuity of $u$, there exist $w,z\in S^{N-1}_{y}$ satisfying
\begin{equation*}\label{yz}
u^*(w)>t\geq u(w)\quad\mathrm{and}\quad u(z)>t\geq u^*(z).
\end{equation*}
Let $H\in\mathcal{H}_0$ with $w\in H$ and $z=\sigma_H(w)$. Since $u^*(w)>u^*(z)$, hence $w$ is closer to $z^1_{y}$ than $z$; that is, $H\in\mathcal{H}_1$. For all $x\in H\cap S^N_+$, using \eqref{invariance}, we have
\begin{eqnarray*}\label{pointhell}
|u_H(x)-u^*(x)|^2&+&|u_H(\sigma_H(x))-u^*(\sigma_H(x))|^2\nonumber\\
&\leq& |u(x)-u^*(x)|^2+|u(\sigma_H(x))-u^*(\sigma_H(x))|^2,
\end{eqnarray*}
and hence also
\begin{eqnarray*}\label{pointhell}
y^a|u_H(x)-u^*(x)|^2&+&y^a|u_H(\sigma_H(x))-u^*(\sigma_H(x))|^2\nonumber\\
&\leq& y^a|u(x)-u^*(x)|^2+y^a|u(\sigma_H(x))-u^*(\sigma_H(x))|^2.
\end{eqnarray*}
By continuity, the inequality is strict in a neighbourhood of $w$. Integrating over $H\cap S^N_+$, \eqref{strictt} follows.
\endproof
For $u\in C(S^N_+)$, the mapping $H\mapsto u_H$ is continuous from $\mathcal{H}_1\sim SO(N)/\mathbb{Z}_2$ to $L^2(S^N_+;\mathrm{d}\mu)$; that is, the polarization depends continuously on its defining half space. A way to see this fact is the following result from \cite{brocksolynin}.
\begin{Lemma}\label{lem2}
Let $u\in C(S^N_+)$ and $\{H_n\}$ be a sequence of half spaces in $\mathcal{H}_1$. If $H\in\mathcal{H}_1$ and
\begin{equation}\label{measco}
\lim_{n\to +\infty}\mu\left((H_n\triangle H)\cap S^N_+\right)=0,
\end{equation}
then $u_{H_n}\longrightarrow u_H$ in $L^2(S^N_+;\mathrm{d}\mu)$.
\end{Lemma}
\proof
By \eqref{measco} we have $\lim_{n\to +\infty}\sigma_{H_n}(z)=\sigma_H(z)$ uniformly on compact subsets of $S^N_+$. Hence the result follows.
\endproof
By compactness of $SO(N)/\mathbb{Z}_2$, if $u\in C(S^N_+)$, the minimization problem
\begin{equation*}\label{minH}
c:=\inf_{H\in\mathcal{H}_1}||u_H-u^*||_{L^2(S^N_+;\mathrm{d}\mu)}
\end{equation*}
is achieved by some $H:=H(u)$.
\begin{Lemma}\label{lem3}
Let $u\in C^{\infty}(S^N_+)$ be a nonnegative function. Then the sequence $\{u_n\}$ defined by $u_0=u$, $u_{n+1}=(u_n)_{H_n}$ and
\begin{equation*}\label{recurH}
||u_{n+1}-u^*||_{L^2(S^N_+;\mathrm{d}\mu)}=\min_{H\in\mathcal{H}_1}||(u_n)_H-u^*||_{L^2(S^N_+;\mathrm{d}\mu)}
\end{equation*}
converges to $u^*$ in $L^2(S^N_+;\mathrm{d}\mu)$.
\end{Lemma}
\proof
Since $u\in C^{\infty}(S^N_+)$ then $u\in W^{1,q}(S^N_+;\mathrm{d}\mu)$ for every $1\leq q<+\infty$ and so for every $n\in\mathbb{N}$ it holds that $||\nabla_{S^N} u_n||_{L^q(S^N_+;\mathrm{d}\mu)}=||\nabla_{S^N} u||_{L^q(S^N_+;\mathrm{d}\mu)}$; that is, the sequence $\{u_n\}$ is bounded in $W^{1,q}(S^N_+;\mathrm{d}\mu)$. Hence, for $q>2$, by the Rellich theorem (compact embedding in H\"older spaces), we can assume, up to a subsequence, that $u_n\longrightarrow v$ uniformly. Since $(u_n)^*=u^*$ and the fact that foliated Schwarz symmetrization is a contraction in $L^p(\mathrm{d}\mu)$-spaces, it follows that $v^*=u^*$. Moreover, for every $H\in\mathcal{H}_1$ we have
\begin{equation}\label{chain1}
||u_{n+1}-u^*||_{L^2(S^N_+;\mathrm{d}\mu)}\leq||(u_n)_H-u^*||_{L^2(S^N_+;\mathrm{d}\mu)}\leq||u_{n}-u^*||_{L^2(S^N_+;\mathrm{d}\mu)},
\end{equation}
where the first inequality follows from our hypothesis and the second one always holds since polarization is a contraction in $L^p(\mathrm{d}\mu)$-spaces. Taking the limit along the subsequence in \eqref{chain1}, we get
\begin{equation*}\label{chain2}
||v-u^*||_{L^2(S^N_+;\mathrm{d}\mu)}\leq||v_H-u^*||_{L^2(S^N_+;\mathrm{d}\mu)}\leq||v-u^*||_{L^2(S^N_+;\mathrm{d}\mu)}.
\end{equation*}
But $v^*=u^*$ and $H\in\mathcal{H}_1$ is arbitrary. So by Lemma \ref{lem1} there are two possibilities: either there exists $H\in\mathcal{H}_1$ such that the second inequality is strict or $v=v^*$. But the first case can't happen and hence the result is proved.
\endproof
As a consequence, we remark that since for every $n\in\mathbb{N}$ the sequence of Lemma \ref{lem3} satisfies $||u_n||_{L^2(S^N_+;\mathrm{d}\mu)}=||u||_{L^2(S^N_+;\mathrm{d}\mu)}$, it holds that
\begin{equation}\label{normu=}
||u^*||_{L^2(S^N_+;\mathrm{d}\mu)}=||u||_{L^2(S^N_+;\mathrm{d}\mu)}.
\end{equation}
Now we can prove the P\'olya-Szeg\"o inequality for the foliated Schwarz symmetrization on the hemisphere.
\begin{Proposition}\label{teo1}
If $u\in H^{1;a}(S^N_+)$ and nonnegative, then $u^*\in H^{1;a}(S^N_+)$, nonnegative, and
\begin{equation}\label{polyaszego}
\int_{S^N_+}y^a|\nabla_{S^N}u^*|^2\leq\int_{S^N_+}y^a|\nabla_{S^N}u|^2.
\end{equation}
\end{Proposition}
\proof
Assume first that $u\in C^{\infty}(S^N_+)$. The sequence $\{u_n\}$ associated to $u$ as in Lemma \ref{lem3} is such that $u_n\longrightarrow u^*$ in $L^2(S^N_+;\mathrm{d}\mu)$ and for every $n\in\mathbb{N}$
\begin{equation*}\label{12}
||u_n||_{L^2(S^N_+;\mathrm{d}\mu)}=||u||_{L^2(S^N_+;\mathrm{d}\mu)}\quad\mathrm{and}\quad||\nabla_{S^N}u_n||_{L^2(S^N_+;\mathrm{d}\mu)}=||\nabla_{S^N}u||_{L^2(S^N_+;\mathrm{d}\mu)}.
\end{equation*}
Hence, $u^*\in H^{1;a}(S^N_+)$ and by the weak lower simicontinuity of the norm in an Hilbert space,  $||\nabla_{S^N}u^*||_{L^2(S^N_+;\mathrm{d}\mu)}\leq||\nabla_{S^N}u||_{L^2(S^N_+;\mathrm{d}\mu)}$.

If $u\in H^{1;a}(S^N_+)$, then by density there exists a sequence $\{u_m\}$ in $C^{\infty}(S^N_+)$ converging to $u$ in $H^{1;a}(S^N_+)$. Since any rearrangement is a contraction in $L^2(\mathrm{d}\mu)$, then $u^*_m\longrightarrow u^*$ in $L^2(S^N_+;\mathrm{d}\mu)$ and hence
\begin{equation*}\label{chain3}
||\nabla_{S^N}u^*||_{L^2(S^N_+;\mathrm{d}\mu)}\leq\liminf_{m\to+\infty}||\nabla_{S^N}u^*_m||_{L^2(S^N_+;\mathrm{d}\mu)}\leq\liminf_{m\to+\infty}||\nabla_{S^N}u_m||_{L^2(S^N_+;\mathrm{d}\mu)}=||\nabla_{S^N}u||_{L^2(S^N_+;\mathrm{d}\mu)}.
\end{equation*}
This completes the proof.
\endproof

\subsection{Optimal geometry for boundary conditions imposing one symmetry}
Let $N=2$ and let us consider $\Sigma_1$ previously defined as a plane containing the axis $y$ with relative reflection $\sigma_1:=\sigma_{\Sigma_1}$ (we remember that we choose the one containing points with angle $\phi=0$). Let us now define the following class of symmetric regions
\begin{equation*}\label{A}
\mathcal{A}_1=\{\omega\subset S^1 : \mathcal{H}^1(\omega)=\mathcal{H}^1(S^1\setminus\omega) \ \mathrm{and} \ (x,0)\in\omega\iff \sigma_{1}(x,0)\in S^1\setminus\omega\}.
\end{equation*}
Hence, we wish to study the problem
\begin{equation}\label{InfInf}
\inf_{\omega\in\mathcal{A}_1}\lambda_1^s(\omega);
\end{equation}
that is, we see the optimal geometry of the boundary condition region $\omega\in\mathcal{A}_1$ as the one which gives the lowest eigenvalue. As we have previously said, for a fixed $\omega\in\mathcal{A}_1$, the minimization of the Rayleigh quotient is standard and we get the existence of a nontrivial and nonnegative minimizer for the energy
\begin{equation*}\label{Ifun}
\int_{S^2_+}y^a|\nabla_{S^2} u|^2
\end{equation*}
constrained to $X_{\omega}=\left\{u\in H^{1;a}_{\omega}(S^2_+) : \int_{S^2_+}y^au^2=1\right\}$. Moreover, the constrained minimizer $u_{\omega}$ found is also a minimizer of the Rayleigh quotient in the whole $H^{1;a}_{\omega}(S^2_+)$. By a simple Frech\'et differentiation of the Rayleigh quotient, turns out to be true that such a minimizer is a weak solution of problem \eqref{Plambdagamma} in the sense that
\begin{equation}\label{weakform}
\int_{S^2_+}y^a\nabla_{S^2}u_{\omega}\nabla_{S^2}\phi=\lambda^s_1(\omega)\int_{S^2_+}y^au_{\omega}\phi,\quad\forall\phi\in C^{\infty}_{0}(S^2_+\cup\omega).
\end{equation}
Thanks to the results obtained for the foliated Schwarz symmetrization, we are able to show the following result.
\begin{Proposition}\label{prop2}
For every fixed $\omega\in\mathcal{A}_1$ let us consider the minimizer $u_{\omega}\in H^{1;a}_{\omega}(S^2_+)$ of the Rayleigh quotient. Then there exists a function $u^*_{\omega}\in H^{1;a}_{\omega_1}(S^2_+)$ such that
\begin{equation*}\label{u*gamma}
\mathcal{R}^a(u^*_{\omega})\leq\mathcal{R}^a(u_{\omega})=\lambda^s_1(\omega),
\end{equation*}
where $\omega_1:=S^1\cap\{0<\phi<\pi\}\in\mathcal{A}_1$ is half of $S^{1}$.
\end{Proposition}
\proof
First we recall that we can choose $u_{\omega}$ nonnegative and it is nontrivial. Then, let us define the function $u^*_{\omega}$ as in Definition \ref{def2}; that is, the foliated Schwarz symmetrization of $u_{\omega}$ so that, on any level $S^1_y$, the decreasing rearrangement is centered in the points $z_y^1$ which has coordinate $\phi=\pi/2$. Hence, thanks to Proposition \ref{teo1}, it holds that
\begin{equation}\label{fili}
\int_{S^2_+}y^a|\nabla_{S^2}u^*_{\omega}|^2\leq\int_{S^2_+}y^a|\nabla_{S^2}u_{\omega}|^2,
\end{equation}
and we know also that
\begin{equation}\label{fili}
\int_{S^2_+}y^a|u^*_{\omega}|^2=\int_{S^2_+}y^a|u_{\omega}|^2;
\end{equation}
that is, the Rayleigh quotient decreases. Moreover, considering the restriction of $u^*_{\omega}$ to $S^1$, we know that the set $\{u^*_{\omega}|_{S^{1}}>0\}$ is the closed geodesic ball centered in $z_0^1$ with measure given by
\begin{equation*}\label{r2}
\mathcal{L}^{1}(\{u^*_{\omega}|_{S^{1}}>0\})=\mathcal{L}^{1}(\{u_{\omega}|_{S^{1}}>0\})=\mathcal{L}^{1}(\omega)=\frac{1}{2}\mathcal{L}^{1}(S^{1}).
\end{equation*}
\endproof
Proposition \ref{prop2} obviously implies that
\begin{equation}\label{infgammatil}
\inf_{\omega\in\mathcal{A}_1}\lambda_1^s(\omega)=\lambda_1^s(\omega_1)=:\lambda_1^s(1),
\end{equation}
and it is attained by a nontrivial and nonnegative minimizer $u_1\in H^{1;a}_{\omega_1}(S^2_+)$ which is a weak solution of
\begin{equation*}\label{Plambda}
\begin{cases}
-L_a^{S^2} u=y^a\lambda_1^s(1) u  &\mathrm{in}\quad S^2_+, \\
u=0  &\mathrm{in}\quad S^{1}\setminus\omega_1, \\
\partial^a_y u=0 &\mathrm{in}\quad\omega_1\subset S^{1},
\end{cases}
\end{equation*}
in the sense of \eqref{weakform}.

\subsection{Optimal geometry for boundary conditions imposing more symmetries}
In this section we wish to show the optimal geometry of the boundary condition region in case of more symmetries; that is, we will consider for an arbitrary $k\in\mathbb{N}$, the boundary condition set $\omega\in\mathcal{A}_k$ where
\begin{equation*}\label{Ak}
\mathcal{A}_k=\{\omega\subset S^{1} : \mathcal{H}^{1}(\omega)=\mathcal{H}^{1}(S^{1}\setminus\omega) \ \mathrm{and} \ (x,0)\in\omega\iff \sigma_{i}(x,0)\in S^1\setminus\omega \ \forall i=1,...,k\},
\end{equation*}
with $\sigma_i:=\sigma_{\Sigma_{i}}$ reflections with respect to $\Sigma_{i}$ planes containing the axis $y$ and hence orthogonal to the plane $y=0$, for every $i=1,...,k$. Considering $T_k=2\pi/k$ as the period, then the plane $\Sigma_{i+1}$ is obtained by rotating $\Sigma_{i}$ with respect to $\phi$ of an angle $T_k/2$.

We are interested in finding solutions $u$ to \eqref{Plambdagamma} with $\omega\in\mathcal{A}_k$ and such that
\begin{equation}\label{syyyy}
u(z)=u(\sigma_{i}(\sigma_{j}(z)))
\end{equation}
for every $i,j=1,...,k$, for almost every $z\in S^2_+$ with respect to the measure given by $\mathrm{d}\mu=y^a\mathrm{d}S(z)$ and also for almost every $z\in S^{1}$ with respect to the $1$-dimensional Lebesgue measure. So, we study the following problem
\begin{equation}\label{InfInfk}
\inf_{\omega\in\mathcal{A}_k}\lambda_1^s(\omega),
\end{equation}
where
\begin{equation}\label{lambdawc}
\lambda_1^s(\omega):=\inf\left\{\mathcal{R}^au \ : \ u\in H^{1;a}_{\omega}(S^2_+)\setminus\{0\} \ \mathrm{and} \ \eqref{syyyy} \ \mathrm{holds}\right\}.
\end{equation}
We remark that the definition of the first eigenvalue with respect to $\omega$ given previously for the case of only one symmetry is in accord with this new definition because \eqref{syyyy} obviously holds in that case.

Let $\omega\in\mathcal{A}_k$. Then there exists a nontrivial and nonnegative minimizer for the functional $\int_{S^2_+}y^a|\nabla u|^2$ constrained to $\overline{X_{\omega}}=\{u\in H^{1;a}_{\omega}(S^2_+) \ : \ \int_{S^2_+}y^au^2=1 \ \mathrm{and} \ \eqref{syyyy} \ \mathrm{holds}\}$. First of all, we remark that the set of functions $\overline{X_{\omega}}$ is not empty. In fact, let us define the fundamental subdomain of $S^2_+$
\begin{equation}\label{sn+t}
S^2_+(k)=\{z\in S^2_+ \ : \ \phi\in(0,T_k)\}.
\end{equation}
Let us now split this domain in other two subdomains $S^2_+(k,1)=\{z\in S^2_+ \ : \ \phi\in(0,T_k/2)\}$ and $S^2_+(k,2)=\{z\in S^2_+ \ : \ \phi\in(T_k/2,T_k)\}$. Since both these domains have positive $L_a$-capacity, we can find two nontrivial nonnegative functions $u_i\in H^{1;a}_0(S^2_+(k,i))$ for $i=1,2$. Then we can merge them in a unique function defined over the fundamental domain and then we can extend it to the whole of $S^2_+$ in a periodic way. After a normalization in $L^2(S^2_+;\mathrm{d}\mu)$, we get an element of $\overline{X_{\omega}}$.

The other thing to remark is that property \eqref{syyyy}, satisfied by the generic minimizing sequence $\{u_n\}_{n=1}^{+\infty}\subseteq \overline{X_{\omega}}$, is also satisfied by its weak limit $u_{\omega}\in H^{1;a}_{\omega}(S^2_+)$, but this fact is trivial using Sobolev embedding in $L^2(S^2_+;\mathrm{d}\mu)$, trace theory in $L^2(S^{1})$, and pointwise convergence. Hence, we wish to show that this critical point $u_{\omega}$ founded minimizing the energy on $\overline{X_{\omega}}$ is also a critical point of the same functional over $X_{\omega}=\{u\in H^{1;a}_{\omega}(S^2_+) \ : \  \int_{S^2_+}y^au^2=1\}$. Let $\mathbb{G}$ be the group of rotation with respect to $\phi$ of a fixed angle $T_k$. Let us consider the action of this group
\begin{eqnarray}\label{group}
&&\mathbb{G}\times X_{\omega}\longrightarrow X_{\omega}\nonumber\\
&&[g,u]\longmapsto u\circ g.
\end{eqnarray}
Since for every $g\in \mathbb{G}$, $g(\omega)=\omega$ and $\int_{S^2_+}y^a|\nabla_{S^2}u\circ g|^2=\int_{S^2_+}y^a|\nabla_{S^2}u|^2$, then the energy is invariant with respect to $\mathbb{G}$ and the action in \eqref{group} is isometric (we remark that the rotation of the group does not change the value in $y$). Hence, by the principle of symmetric criticality of Palais, a critical point of the energy over the set
\begin{equation*}\label{fix}
Fix(\mathbb{G})=\{u\in X_{\omega} \ : \ u\circ g=u \ \forall g\in \mathbb{G}\}=\overline{X_{\omega}},
\end{equation*}
is also a critical point of the same functional over $X_{\omega}$. Then, it follows easily that $u_{\omega}$ is also a critical point of the Rayleigh quotient over the whole $H^{1;a}_{\omega}(S^2_+)$; that is, it is a solution to \eqref{Plambdagamma} with $\omega\in\mathcal{A}_k$ and such that property \eqref{syyyy} holds.

Hence, by the symmetry condition \eqref{syyyy}, if we know $u_{\omega}$ in $S^2_+(k)$, then $u_{\omega}$ is consequently determined in the whole hemisphere $S^2_+$. To simplify the notation let us call $u:=u_{\omega}$. Let us define over the whole hemisphere the function
\begin{equation}\label{v}
v(\theta,\phi):=u(\theta,\phi/k).
\end{equation}
Obviously $v\in H^{1;a}_{\overline{\omega}}(S^2_+)$ with $\overline{\omega}\in\mathcal{A}_1$ and it is nonnegative. Following the same steps done before, we wish to rearrange the function $v$, in order to lower the $L^2(\mathrm{d}\mu)$-norm of its tangential gradient, by the foliated Schwarz hemispherical symmetrization. Actually we will consider a gradient-type operator such that
\begin{equation}\label{gradk}
|\nabla^{(k)}_{S^2}v|^2:=(\partial_{\theta}v)^2+\frac{k^2}{y^2}(\partial_{\phi}v)^2.
\end{equation}
The following P\'olya-Szeg\"o type inequality holds.
\begin{Proposition}\label{teo2}
Let us consider $v^*$ as the foliated Schwarz symmetrization of the function $v\in H^{1;a}_{\overline{\omega}}(S^2_+)$ defined in \eqref{v}. Then $v^*\in H^{1;a}(S^2_+)$ and
\begin{equation*}\label{psv}
\int_{S^2_+}y^a|\nabla^{(k)}_{S^2}v^*|^2\leq\int_{S^2_+}y^a|\nabla^{(k)}_{S^2}v|^2.
\end{equation*}
\end{Proposition}
\proof
Following the same steps seen in Lemma \ref{lem3} for the case $k=1$, if $v\in C^{\infty}(S^2_+)$, then we construct the sequence $\{v_n\}$ of polarized functions such that $v_n\longrightarrow v^*$ in $L^2(S^2_+;\mathrm{d}\mu)$, where $v^*$ is defined as in the proof of Proposition \ref{prop2}. In \cite{hajaiejsquass} it is proved that for every $p\in(1,+\infty)$ and for every suitable half space, one has
\begin{equation*}\label{partialH}
||D_iv||_{L^p(S^2_+)}=||D_iv_H||_{L^p(S^2_+)},
\end{equation*}
for every first order derivative; that is,
\begin{equation}\label{parttt}
||\partial_{\theta}v||_{L^2(S^2_+)}=||\partial_{\theta}v_H||_{L^2(S^2_+)}\quad\mathrm{and}\quad||\partial_{\phi}v||_{L^2(S^2_+)}=||\partial_{\phi}v_H||_{L^2(S^2_+)}.
\end{equation}
From \eqref{parttt}, it follows that also
\begin{equation}\label{teta}
\int_{S^2_+}y^{a-2}(\partial_{\phi}v)^2=\int_{S^2_+}y^{a-2}(\partial_{\phi}v_H)^2,
\end{equation}
since it holds that for every point $z\in S^2_+$, the point $\sigma_H(z)$ has the same coordinate $y$. Then, by \eqref{parttt} and \eqref{teta} it follows that
\begin{equation*}\label{gradk=H}
\int_{S^2_+}y^a|\nabla^{(k)}_{S^2}v_H|^2=\int_{S^2_+}y^a|\nabla^{(k)}_{S^2}v|^2.
\end{equation*}
Moreover, it is easy to see that the quantity $\int_{S^2_+}y^a|\nabla^{(k)}_{S^2}v|^2$ is an equivalent norm on $H^{1;a}_{\overline{\omega}}(S^2_+)$; that is,
\begin{equation*}\label{equnorm}
\int_{S^2_+}y^a|\nabla_{S^2}v|^2\leq\int_{S^2_+}y^a|\nabla^{(k)}_{S^2}v|^2\leq k^2\int_{S^2_+}y^a|\nabla_{S^2}v|^2.
\end{equation*}
Hence, using the weak lower semicontinuity of the norm on an Hilbert space, we can replicate the proof of Proposition \ref{teo1} using the new gradient-type norm. Working first with $v\in C^{\infty}(S^2_+)$ and then in $H^{1;a}_{\overline{\omega}}(S^2_+)$ by a density argument, the result is easily proved.
\endproof
Since
\begin{eqnarray*}\label{z1}
|\nabla^{(k)}_{S^2}v(\theta,\phi)|^2&=&\left(\partial_{\theta}[u(\theta,\phi/k)]\right)^2+\frac{k^2}{y^2}\left(\partial_{\phi}[u(\theta,\phi/k)]\right)^2\nonumber\\
&=&\left(u_{\theta}(\theta,\phi/k)\right)^2+\frac{k^2}{y^2}\left(\frac{1}{k}u_{\phi}(\theta,\phi/k)\right)^2\nonumber\\
&=&|\nabla_{S^2}u(\theta,\phi/k)|^2,
\end{eqnarray*}
hence it holds that
\begin{equation*}\label{psk}
\int_{S^2_+}y^a|\nabla_{S^2}u^*(\theta,\phi/k)|^2\leq\int_{S^2_+}y^a|\nabla_{S^2}u(\theta,\phi/k)|^2,
\end{equation*}
and changing variables we get that
\begin{equation*}\label{psP}
\int_{S^2_+(k)}ky^a|\nabla_{S^2}u^*|^2\leq\int_{S^2_+(k)}ky^a|\nabla_{S^2}u|^2.
\end{equation*}
Obviously $u^*$ defines a unique function, thanks to condition \eqref{syyyy}, over $S^2_+$ and it is easy to check that
\begin{equation*}\label{kvolte}
\int_{S^2_+}y^a|\nabla_{S^2}u|^2=\int_{S^2_+(k)}ky^a|\nabla_{S^2}u|^2\quad\mathrm{and}\quad\int_{S^2_+}y^a|\nabla_{S^2}u^*|^2=\int_{S^2_+(k)}ky^a|\nabla_{S^2}u^*|^2.
\end{equation*}
Moreover, this fact says us that $u^*\in H^{1;a}_{\omega_k}(S^2_+)$ where $\omega_k:=S^1\cap\{\phi\in\bigcup_{i=1}^{k}((i-1)T_k,(i-1/2)T_k)\}\in\mathcal{A}_k$ is the particular boundary condition set that is the most connected one, according with the conditions given. Finally it follows easily that $R^a(u^*)\leq R^a(u_{\omega})=\lambda^s_1(\omega)$; that is,
\begin{equation}\label{infgammatilk}
\inf_{\omega\in\mathcal{A}_k}\lambda^s_1(\omega)=\lambda^s_1(\omega_k)=:\lambda^s_1(k),
\end{equation}
in the sense of \eqref{lambdawc}. Moreover, the minimization problem in \eqref{infgammatilk} admits a nontrivial and nonnegative minimizer $u_k\in H^{1;a}_{\omega_k}(S^2_+)$, which is also a weak solution of
\begin{equation*}\label{Plambdak}
\begin{cases}
-L_a^{S^2} u=y^a\lambda_1^s(k) u  &\mathrm{in}\quad S^2_+, \\
u=0  &\mathrm{in}\quad S^{1}\setminus\omega_k, \\
\partial^a_y u=0 &\mathrm{in}\quad\omega_k\subset S^{1},
\end{cases}
\end{equation*}
in the sense of \eqref{weakform} and such that condition \eqref{syyyy} is satisfied.

\subsection{Ordering eigenvalues with respect to the number of symmetries}
The aim of this section is to show that the sequence of eigenvalues $\{\lambda^s_1(k)\}_{k=1}^{+\infty}$, obtained for every $k\in\mathbb{N}$ optimizing the energy under the best boundary condition, is such that
\begin{equation}\label{chainlambda}
\lambda_1^s(S^1)\leq\lambda_1^s(1)\leq...\leq\lambda_1^s(k)\leq\lambda_1^s(k+1)\leq...\leq\lambda_1^s(\emptyset).
\end{equation}
First, we remark that by \eqref{zeig}, then for every $k\in\mathbb{N}$ it holds that $$\lambda_1^s(S^1)\leq\lambda_1^s(k)\leq\lambda_1^s(\emptyset).$$
Let $k\in\mathbb{N}$ fixed and $\omega\in\mathcal{A}_k$. Let us define $u=u_{\omega}$ the minimizer for the problem \eqref{lambdawc} and $v$ as in \eqref{v}. Then, we have proved that
\begin{equation*}\label{chain4}
\int_{S^2_+}y^a|\nabla^{(k)}_{S^2}v|^2=\int_{S^2_+}y^a|\nabla_{S^2}u(\theta,\phi/k)|^2=\int_{S^2_+(k)}ky^a|\nabla_{S^2}u|^2=\int_{S^2_+}y^a|\nabla_{S^2}u|^2.
\end{equation*}
Hence, the eigenvalue $\lambda_1^s(k)$ can be also expressed as
\begin{equation*}\label{infff}
\lambda_1^s(k)=\inf_{\omega\in\mathcal{A}_1}\left(\inf\left\{\int_{S^2_+}y^a|\nabla^{(k)}_{S^2}u|^2 \ : \ u\in H^{1;a}_{\omega}(S^2_+) \ \mathrm{with} \ \int_{S^2_+}y^au^2=1\right\}\right),
\end{equation*}
and this quantity is obviously non decreasing in $k\in\mathbb{N}$. This implies \eqref{chainlambda}.

From now on, let us consider the sequence $\{u_k\}_{k=1}^{+\infty}\subseteq H^{1;a}(S^2_+)$ of nonnegative first eigenfunctions associated to the sequence $\{\lambda_1^s(k)\}_{k=1}^{+\infty}$ and such that
\begin{equation}\label{uk}
\int_{S^2_+}y^a|\nabla_{S^2}u_k|^2=\lambda_1^s(k)\quad\mathrm{and}\quad\int_{S^2_+}y^au_k^2=1/2.
\end{equation}

\subsection{H\"older regularity of eigenfunctions}
We remark that the minimization problem under $k$ symmetries seen in \eqref{infgammatilk} can be extended in a natural way, in the case of two components which are segregated on $S^1$ and satisfy some symmetry and measure conditions. Let us define the set of 2-partitions of $S^1$ satisfying a condition over the measure and one over the symmetry
\begin{eqnarray}\label{part}
\mathcal{P}^2_k &=&\{(\omega_1,\omega_2) : \ \omega_i\subset S^{1} \mathrm{ \ open}, \ \omega_1\cap\omega_2=\emptyset,\nonumber\\
&&\overline{\omega_1}\cup\overline{\omega_2}= S^{1}, \ \mathcal{H}^1(\omega_1)=\mathcal{H}^1(\omega_2), \ z\in\omega_1\Leftrightarrow \sigma_i(z)\in\omega_2 \ \forall i=1,...,k\}.
\end{eqnarray}
Fixing a couple $(\omega_1,\omega_2)\in\mathcal{P}^2_k$, let us also define the set of functions
\begin{eqnarray}\label{Bk}
\mathcal{B}_k(\omega_1,\omega_2)&=&\{(u_1,u_2): \ u_i\in H^{1;a}(S^2_+), \ \int_{S^2_+}y^au_i^2=1, \ u_i=0 \ \mathrm{in \ } S^1\setminus \omega_i, \ \mathrm{with \ } (\omega_1,\omega_2)\in\mathcal{P}^2_k,\nonumber\\ 
&& u_i(z)=u_i(\sigma_j(\sigma_l(z))) \ \mathrm{and} \ u_1(z)=u_2(\sigma_j(z)) \ \mathrm{in \ } S^2_+,\nonumber\\
 &&\forall i=1,2, \ j,l=1,...,k\}.
\end{eqnarray}
First of all, we remark that also in this case it is easy to check that, for any fixed couple $(\omega_1,\omega_2)\in\mathcal{P}^2_k$, the set $\mathcal{B}_k(\omega_1,\omega_2)$ is not empty. In fact, proceeding as in section 3.3, we first construct the first component $u_1$ on the fundamental domain $S^2_+(k)$ and then we extend it in a periodic way over $S^2_+$ and we normalize it in $L^2(S^2_+;\mathrm{d}\mu)$. Hence, we can define the second component $u_2$ such that $u_2(z)=u_1(\sigma_i(z))$ for any $i=1,...,k$.

So, as it happened in \eqref{infgammatilk} for the case of one component, we consider the minimization problem
\begin{equation}\label{sumlambda}
\inf_{(\omega_1,\omega_2)\in\mathcal{P}^2_k}\inf_{(u_1,u_2)\in\mathcal{B}_k(\omega_1,\omega_2)}I(u_1,u_2),
\end{equation}
where
\begin{equation}
I(u_1,u_2)=\frac{1}{2}\int_{S^2_+}y^a\left(|\nabla_{S^2}u_1|^2+|\nabla_{S^2}u_2|^2\right).
\end{equation}
Hence, the problem in \eqref{sumlambda} is equivalent to
\begin{equation}
\inf_{(\omega_1,\omega_2)\in\mathcal{P}^2_k}\frac{\lambda_1^s(\omega_1)+\lambda_1^s(\omega_2)}{2}.
\end{equation}
Working with the foliated Schwarz symmetrization on both the components, with respect to both the opposite poles $z_0^1$ and $z_0^2$, it happens that the infimum is achieved by the couple $(u_k,v_k)$ where $u_k$ is the minimizer of $\lambda_1^s(\omega_k)$ found for the problem \eqref{infgammatilk}, $\omega_k:=S^1\cap\{\phi\in\bigcup_{i=1}^{k}((i-1)T_k,(i-1/2)T_k)\}\in\mathcal{A}_k$, and $v_k$ is such that $v_k(z)=u_k(\sigma_j(z))$ in $S^2_+$ for every $j=1,...,k$; that is, $v_k$ achieves $\lambda_1^s(\omega_k^c)$. Moreover, the infimum in \eqref{sumlambda} is given by the number 
\begin{equation}
\inf_{(\omega_1,\omega_2)\in\mathcal{P}^2_k}\frac{\lambda_1^s(\omega_1)+\lambda_1^s(\omega_2)}{2}=\frac{\lambda_1^s(\omega_k)+\lambda_1^s(\omega_k^c)}{2}=\lambda_1^s(k).
\end{equation}
Let us define 
\begin{eqnarray}
X &=&\{(u_1,u_2): \ u_i\in H^{1;a}(S^2_+), \ \int_{S^2_+}y^au_i^2=1, \ u_1=0 \ \mathrm{in \ } S^1\setminus \omega_k,\nonumber\\
 && u_2=0 \ \mathrm{in \ } S^1\setminus \omega_k^c,\mathrm{with \ } (\omega_k,\omega_k^c)\in\mathcal{P}^2_k\},
\end{eqnarray}
and also the group $\mathbb{G}$ of all the reflections $\sigma_i$, with $i=1,...,k$ endowed with the composition between reflections. Let us define the action
\begin{eqnarray}\label{group3}
&& X\times \mathbb{G}\longrightarrow X\nonumber\\
&&[(u_1,u_2),g]\longmapsto (u_2\circ g,u_1\circ g).
\end{eqnarray}
That is, for $g=\sigma_i$, it holds $$[(u_1,u_2),\sigma_i]=(u_2\circ \sigma_i,u_1\circ \sigma_i),$$ and for $g=\sigma_i\circ\sigma_j$, it holds $$[(u_1,u_2),\sigma_i\circ\sigma_j]=[[(u_1,u_2),\sigma_i],\sigma_j]=[(u_2\circ \sigma_i,u_1\circ \sigma_i),\sigma_j]=(u_1\circ \sigma_i\circ\sigma_j,u_2\circ \sigma_i\circ\sigma_j).$$
It is easy to check that this action is isometric and that the functional $I(u_1,u_2)$ is invariant with respect to this action. Since $\mathcal{B}_k(\omega_k,\omega_k^c)=Fix(\mathbb{G})$, by the principle of symmetric criticality of Palais, the minimizer $(u_k,v_k)$ is also a nonnegative critical point for $I$ over the whole $X$ and hence a weak solution to the problem
\begin{equation}\label{asd}
\begin{cases}
-L_a^{S^2} u_k=y^a\lambda_1^s(k) u_k, \ -L_a^{S^2} v_k=y^a\lambda_1^s(k) v_k & \mathrm{in}\quad S^2_+, \\
u_k\partial^a_y u_k=0, \ v_k\partial^a_y v_k=0 & \mathrm{in}\quad S^1, \\
u_kv_k=0, & \mathrm{in}\quad S^1.
\end{cases}
\end{equation}
We wish to prove the $C^{0,\alpha}(\overline{S^2_+})$-regularity for $(u_k,v_k)$ via the convergence of solutions of $\beta$-problems over $S^2_+$ to our eigenfunctions.
Let us now consider the following set of functions
\begin{eqnarray}\label{Ck}
\mathcal{C}_k&=&\{(u_1,u_2): \ u_i\in H^{1;a}(S^2_+), \ \int_{S^2_+}y^au_i^2=1, \ u_i(z)=u_i(\sigma_j(\sigma_l(z))) \nonumber\\ 
&&\mathrm{and} \ u_1(z)=u_2(\sigma_j(z)) \ \mathrm{in \ } S^2_+, \ \forall i=1,2, \ j,l=1,...,k\}.
\end{eqnarray}
This space is trivially not empty since $(\mu(S^2_+)^{-1},\mu(S^2_+)^{-1})\in\mathcal{C}_k$.

Hence, for any $\beta>0$, we consider the following minimizization problem
\begin{equation}\label{minbetalambda}
\inf_{(u_1,u_2)\in\mathcal{C}_k}J_\beta(u_1,u_2),
\end{equation}
with
\begin{equation}
J_\beta(u_1,u_2)=\frac{1}{2}\int_{S^2_+}y^a\left(|\nabla_{S^2}u_1|^2+|\nabla_{S^2}u_2|^2\right)+\frac{1}{2}\int_{S^1}\beta u_1^2u_2^2=I(u_1,u_2)+\frac{1}{2}\int_{S^1}\beta u_1^2u_2^2.
\end{equation}
For every $\beta>0$ fix, the functional $J_\beta$ is Gateaux derivable in any direction, coercive and weakly lower semicontinuous in $\mathcal{C}_k$, and hence there exists a nonnegative minimizer $(u_\beta,v_\beta)\in\mathcal{C}_k$. Moreover by the previous argument, defining
\begin{equation}
Y=\left\{(u_1,u_2): \ u_i\in H^{1;a}(S^2_+), \ \int_{S^2_+}y^au_i^2=1\right\},
\end{equation}
since $J_\beta$ is invariant with respect to the action $Y\times \mathbb{G}\longrightarrow Y$ with $\mathbb{G}$ as in \eqref{group3}, we get that this minimizer is also a critical point over $Y$ and hence a weak solution to
\begin{equation}\label{asd2}
\begin{cases}
-L_a^{S^2} u_\beta=y^a\lambda_\beta u_\beta, \ -L_a^{S^2} v_\beta=y^a\lambda_\beta v_\beta & \mathrm{in}\quad S^2_+, \\
-\partial^a_y u_\beta=\beta u_\beta v^2_\beta, \ -\partial^a_y v_\beta=\beta v_\beta u^2_\beta & \mathrm{in}\quad S^1,
\end{cases}
\end{equation}
where $\lambda_\beta=\int_{S^2_+}y^a|\nabla_{S^2}u_\beta|^2+\int_{S^1}\beta u_\beta^2v_\beta^2=\int_{S^2_+}y^a|\nabla_{S^2}v_\beta|^2+\int_{S^1}\beta u_\beta^2v_\beta^2$.
Moreover, since the couple $(u_k,v_k)\in\mathcal{B}_k(\omega_k,\omega_k^c)\subset\mathcal{C}_k$, it holds that for any $\beta>0$, we get the uniform bound
\begin{equation}\label{bou}
0\leq\lambda_\beta\leq 2J_\beta(u_\beta,v_\beta)\leq 2J_\beta(u_k,v_k)=2\lambda_1^s(k).
\end{equation}
This uniform bound gives the weak convergence in $H^{1;a}(S^2_+)$ of the $\beta$-sequence to a function $(u_\infty,v_\infty)$. Moreover, since solutions to \eqref{asd2} are bounded in $C^{0,\alpha}(\overline{S^2_+})$ uniformly in $\beta>0$ for $\alpha>0$ small, as it is proved in \cite{terverzil2}, we obtain, up to consider a subsequence as $\beta\to+\infty$, that the convergence is uniform on compact sets and so that the limit satisfies the symmetry conditions. Moreover it holds that
\begin{equation}\label{bou2}
0\leq\lambda_\beta= J_\beta(u_\beta,v_\beta)+\frac{1}{2}\int_{S^1}\beta u_\beta^2v_\beta^2\leq\lambda_1^s(k)+\frac{1}{2}\int_{S^1}\beta u_\beta^2v_\beta^2,
\end{equation}
and since $\frac{1}{2}\int_{S^1}\beta u_\beta^2v_\beta^2\to 0$ (see Lemma 4.6 in \cite{terverzil2} and Lemma 5.6 in \cite{terverzil} for the details in the case $s=1/2$), the limit should have the two components segregated on $S^1$; that is, $(u_\infty,v_\infty)\in\mathcal{B}_k(\omega_k,\omega_k^c)$ (by the symmetries), and by the minimality of $(u_k,v_k)$ and \eqref{bou2}, we obtain that $(u_\infty,v_\infty)$ owns the same norm of $(u_k,v_k)$ in $H^{1;a}(S^2_+)$, and hence we can choose as a minimizer $(u_\infty,v_\infty)$ which inherits the H\"older regularity up to the boundary.

\subsection{The limit for $k\to+\infty$}
Hence, we have found for any $k\in\mathbb{N}$ fix, a couple $(u_k,v_k)$ of nonnegative eigenfunctions related to $\lambda_1^s(k)$ with the desired symmetry properties. Moreover, for these eigenfunctions we have the regularity $C^{0,\alpha}(\overline{S^2_+})$. Then, we will study the convergence of the sequence of normalized eigenfunctions associated to $\{\lambda_1^s(k)\}_{k=1}^{+\infty}$.

By \eqref{chainlambda} and \eqref{uk}, the sequence $\{u_k\}_{k=1}^{+\infty}$ is uniformly bounded in $H^{1;a}(S^2_+)$ and hence we get, up to consider a subsequence, weak convergence to a function $u$ in $H^{1;a}(S^2_+)$, strong convergence in $L^2(S^2_+;\mathrm{d}\mu)$ with $\int_{S^2_+}y^au^2=1/2$ (we can always renormalize $\{u_k\}_{k=1}^{+\infty}$ so that $\int_{S^2_+}y^au_k^2=1/2$), and pointwise convergence in $S^2_+$ almost everywhere with respect to $\mu$. Moreover, by trace theory we have $L^2(S^1)$-strong convergence on the boundary $S^1$ and also pointwise convergence almost everywhere in $S^1$ with respect to the $1$-dimensional Lebesgue measure. For every $\varepsilon>0$ it holds that $|u(x)|<\varepsilon$ for almost every $x\in S^{1}$ with respect to the $1$-dimensional Lebesgue measure; that is, $u=0$ in $S^{1}$. In fact, fixed $\varepsilon>0$ and $x\in S^{1}$, there exists a $k\in\mathbb{N}$ big enough such that
\begin{equation*}\label{small1}
|u(x)-u_k(x)|<\varepsilon
\end{equation*}
by the pointwise convergence in $S^1$, and such that
\begin{equation}\label{small2}
M|x-\sigma_{\overline{i}}(x)|^{\alpha}<\varepsilon,
\end{equation}
where $M>0$ is a constant, $\alpha$ is the H\"older continuity exponent and $\sigma_{\overline{i}}(x)\in S^{1}$ is the reflection of the point $x$ with respect to the closest symmetrizing plane $\Sigma_{\overline{i}}$. Obviously \eqref{small2} holds because for a $k\in\mathbb{N}$ big enough we can make the distance $|x-\sigma_{\overline{i}}(x)|$ arbitrarily small. Moreover $u_k(\sigma_{\overline{i}}(x))=0$. Hence,
\begin{eqnarray*}\label{u=0}
|u(x)|&=&|u(x)-u_k(x)+u_k(x)-u_k(\sigma_{\overline{i}}(x))|\nonumber\\
&\leq&|u(x)-u_k(x)|+|u_k(x)-u_k(\sigma_{\overline{i}}(x))|\nonumber\\
&\leq&|u(x)-u_k(x)|+M|x-\sigma_{\overline{i}}(x)|^{\alpha}\nonumber\\
&<& 2\varepsilon.
\end{eqnarray*}
Now, we wish to prove that the limit $u$ is a first nonnegative and nontrivial eigenfunction related to $\lambda_1^s(\emptyset)$. First, by the weak convergence of $u_k$ to $u$ in $H^{1;a}(S^2_+)$ and the fact that the limit is such that $u=0$ in $S^{1}$, we get that $u\in H^{1;a}_{0}(S^2_+)$. Moreover, since $C^{\infty}_0(S^2_+)\subseteq C^{\infty}_0(S^2_+\cup\omega_k)$ for every $k\in\mathbb{N}$ and fixing $k\in\mathbb{N}$ it holds that
\begin{equation*}\label{weakformk}
\int_{S^2_+}y^a\nabla_{S^2}u_k\nabla_{S^2}\phi=\lambda_1^s(k)\int_{S^2_+}y^au_k\phi,\quad\forall\phi\in C^{\infty}_0(S^2_+\cup\omega_k),
\end{equation*}
obviously for every $k\in\mathbb{N}$ we obtain that
\begin{equation}\label{weakkkk}
\int_{S^2_+}y^a\nabla_{S^2}u_k\nabla_{S^2}\phi=\lambda_1^s(k)\int_{S^2_+}y^au_k\phi,\quad\forall\phi\in C^{\infty}_0(S^2_+).
\end{equation}
Since the sequence $\{\lambda_1^s(k)\}_{k=1}^{+\infty}$ is non decreasing and bounded from above by $\lambda_1^s(\emptyset)>0$, then
\begin{equation}\label{limlam}
\lim_{k\to+\infty}\lambda_1^s(k)=\tilde{\lambda}\leq\lambda_1^s(\emptyset).
\end{equation}
The weak convergence in $H^{1;a}(S^2_+)$ means that
\begin{equation}\label{weakconv2}
\int_{S^2_+}y^au_k\phi \ + \ \int_{S^2_+}y^a\nabla_{S^2}u_k\nabla_{S^2}\phi\longrightarrow\int_{S^2_+}y^au\phi \ + \ \int_{S^2_+}y^a\nabla_{S^2}u\nabla_{S^2}\phi\quad\forall\phi\in H^{1;a}(S^2_+).
\end{equation}
Since, up to a subsequence, $u_k\longrightarrow u$ in $L^2(S^2_+;\mathrm{d}\mu)$, then it holds also that $u_k\rightharpoonup u$ in $L^2(S^2_+;\mathrm{d}\mu)$; that is,
\begin{equation}\label{weakconv3}
\int_{S^2_+}y^au_k\phi\longrightarrow\int_{S^2_+}y^au\phi\quad\forall\phi\in L^2(S^2_+;\mathrm{d}\mu).
\end{equation}
Since $C^{\infty}_0(S^2_+)\subseteq H^{1;a}(S^2_+)\subseteq L^2(S^2_+;\mathrm{d}\mu)$, then obviously \eqref{weakconv2} and \eqref{weakconv3} hold for every $\phi\in C^{\infty}_0(S^2_+)$. Finally, passing to the limit for $k$ that goes to infinity in \eqref{weakkkk} and putting together \eqref{limlam}, \eqref{weakconv2} and \eqref{weakconv3}, it happens that $u\in H^{1;a}_0(S^2_+)$ satisfies
\begin{equation*}\label{weakkkku}
\int_{S^2_+}y^a\nabla_{S^2}u\nabla_{S^2}\phi=\tilde{\lambda}\int_{S^2_+}y^au\phi,\quad\forall\phi\in C^{\infty}_0(S^2_+);
\end{equation*}
that is, $u$ is an eigenfunction of the problem \eqref{Plambdagamma} with boundary condition $\omega=\emptyset$. Hence $\tilde{\lambda}$ is an eigenvalue of this problem with $\tilde{\lambda}\geq\lambda_1^s(\emptyset)$ since $\lambda_1^s(\emptyset)$ is by definition the smallest one with this boundary condition. Then, by \eqref{limlam}, we get that $\tilde{\lambda}=\lambda_1^s(\emptyset)$.

\subsection{Existence of solutions on the unit half ball}
Our aim is to construct some positive solutions to \eqref{cafsilvess} in case $N=2$ related with the symmetries imposed for the hemispherical problem \eqref{Plambdagamma}. Such solutions will have asymptotic growth rate at infinity which is arbitrarily close to the critical one; that is, $2s$.

Since we have gained H\"older regularity, by \eqref{coloper}, we remark that the first and the last inequalities in the chain \eqref{chainlambda} are strict. In fact, for any $k\in\mathbb{N}$ it holds $\emptyset\subset\omega_k\subset S^1$ and $\mathcal{H}^1(S^1)>\mathcal{H}^1(\omega_k)>\mathcal{H}^1(\emptyset)=0$, and hence
\begin{equation}\label{chainlambda2}
\lambda_1^s(S^1)<\lambda_1^s(1)\leq...\leq\lambda_1^s(k)\leq\lambda_1^s(k+1)\leq...<\lambda_1^s(\emptyset).
\end{equation}
Let us define for every fixed number of symmetries $k\in\mathbb{N}$ the characteristic exponent
\begin{equation}\label{alphak}
d(k):=d(\lambda_1^s(k))=\sqrt{\left(\frac{N-2s}{2}\right)^2+\lambda_1^s(k)}-\frac{N-2s}{2},
\end{equation}
where the sequence of first eigenvalues $\{\lambda_1^s(k)\}$ is defined in section 3.2 and 3.3. Obviously by \eqref{chainlambda} it follows that the degree $d(k)$ is non decreasing in $k$ and in \cite{terverzil2} it is proved that $d(1)=s$. Hence,
\begin{equation}\label{chainalpha}
s=d(1)\leq...\leq d(k)\leq d(k+1)\leq...<d(\lambda_1^s(\emptyset))=2s.
\end{equation}
Therefore, by the previous section, we know that $d(k)\longrightarrow 2s$ as $k\to+\infty$.

From now on, we will follow some ideas and constructions contained in \cite{bereterra,soazil1} for the local case. Now, for every fixed $k\in\mathbb{N}$ and $\beta>1$, we wish to construct over $B_1^+\subset\mathbb{R}^3_+$ nonnegative solutions to
\begin{equation}\label{systb1}
\begin{cases}
L_a u=L_a v=0 & \mathrm{in}\quad B^+_1, \\
-\partial^a_y u=\beta uv^2, \ -\partial^a_y v=\beta vu^2 & \mathrm{in}\quad \partial^0B^+_1, \\
u=g_k, \ v=h_k & \mathrm{in}\quad \partial^+B^+_1,
\end{cases}
\end{equation}
where $(g_k,h_k)\in\mathcal{B}_k$ are nonnegative nontrivial eigenfunctions related to $\lambda_1^s(k)$ satisfying \eqref{asd} and hence such that it holds
\begin{equation}\label{syyyya}
g_k(z)=h_k(\sigma_{i}(z))
\end{equation}
for every $i=1,...,k$. Moreover we choose eigenfunctions as in \eqref{uk} and hence with the property
\begin{equation}\label{h=1}
\int_{\partial^+B^+_1}y^a(g_k^2+h_k^2)=1.
\end{equation}
For simplicity of notations, from now on let us redefine $\lambda=\lambda_1^s(k)$, $d=d(k)$, $g=g_k$, $h=h_k$ and as before $\sigma_i=\sigma_{\Sigma_i}$ the reflection with respect to plane $\Sigma_i$ for every $i=1,...,k$.
\begin{Lemma}\label{teo3}
There exists a pair of nonnegative solutions $(u_\beta,v_\beta)$ to problem \eqref{systb1} satisfying
\begin{enumerate}
\item[1.] for every $i,j=1,...,k$
\begin{equation}\label{symmmm}
\left\{
\begin{aligned}
&u_\beta(z)=u_\beta(\sigma_i(\sigma_j(z))), \\
&v_\beta(z)=v_\beta(\sigma_i(\sigma_j(z))), \\
&u_\beta(z)=v_\beta(\sigma_i(z)); \\
\end{aligned}
\right.
\end{equation}
\item[2.] letting
\begin{equation}\label{functball}
I(u,v):=\frac{1}{2}\int_{B_1^+}y^a(|\nabla u|^2+|\nabla v|^2) \ + \frac{1}{2}\ \int_{\partial^0B_1^+}\beta u^2v^2,
\end{equation}
the uniform estimate $2I(u_\beta,v_\beta)\leq d$ holds.
\end{enumerate} 
\end{Lemma}
\proof
First of all, let us consider in $B^+_1$ the functions
\begin{equation}\label{omoalpk}
\left(G(z),H(z)\right):=|z|^{d}\left(g\left(\frac{z}{|z|}\right),h\left(\frac{z}{|z|}\right)\right),
\end{equation}
which are the $d$-homogeneous extension of $(g,h)$. Since $g,h\in H^{1;a}(S^2_+)$, then it follows by simple calculations that $G,H\in H^{1;a}(B^+_1)$. A weak solution to \eqref{systb1} has to satisfy the following weak formulation
\begin{equation}\label{weak+r}
\left\{
\begin{aligned}
&\int_{B_1^+}y^a\nabla u\nabla \phi \ + \ \int_{\partial^0B_1^+}\beta uv^2\phi=0, \\
&\int_{B_1^+}y^a\nabla v\nabla \phi \ + \ \int_{\partial^0B_1^+}\beta vu^2\phi=0, \\
\end{aligned}
\right.
\end{equation}
for every $\phi\in H^{1;a}_{\partial^+B^+_1}(B^+_1):=\{u\in H^{1;a}(B^+_1) \ : \ u=0 \ \mathrm{in} \ \partial^+B^+_1\}$. Hence, a weak solution to \eqref{systb1} is also a critical point of the functional defined in \eqref{functball} over the reflexive Banach space
\begin{equation}\label{X}
X:=\left(G+H^{1;a}_{\partial^+B^+_1}(B^+_1)\right)\times\left(H+H^{1;a}_{\partial^+B^+_1}(B^+_1)\right).
\end{equation}
In order to get condition 1, we  wish to minimize $I$ over a closed subspace of $X$; that is, $\mathcal{U}\subset X$ the set of pairs of nonnegative functions $(u,v)$ satisfying condition 1. Proceeding as in section 3.5 it is easy to see that $\mathcal{U}$ is not empty. Obviously also $\mathcal{U}$ is a reflexive Banach space and hence, by the direct method of the Calculus of Variations, we have only to show that $I$ is G\^ateaux differentiable in any direction $\phi\in H^{1;a}_{\partial^+B^+_1}(B^+_1)$ such that $(\phi+G,\phi+H)\in\mathcal{U}$, coercive and weakly lower semicontinuous, in order to find a minimizer. The differentiability is a standard calculation that gives us the desired condition
\begin{equation}\label{deri}
\frac{\partial I}{\partial u}(u,v)[\phi]=\int_{B_1^+}y^a\nabla u\nabla \phi \ + \ \int_{\partial^0B_1^+}\beta uv^2\phi\qquad\mathrm{and}\qquad\frac{\partial I}{\partial v}(u,v)[\phi]=\int_{B_1^+}y^a\nabla v\nabla \phi \ + \ \int_{\partial^0B_1^+}\beta vu^2\phi,
\end{equation}
for every direction $\phi\in H^{1;a}_{\partial^+B^+_1}(B^+_1)$ such that $(\phi+G,\phi+H)\in\mathcal{U}$.

Let us recall that $\mathcal{U}$, as a closed subspace, inherits the topology from $X$; that is, the convergence of a pair is characterized by the convergence of its components. Hence, the weak convergence $(u_n,v_n)\rightharpoonup (u,v)$ in $\mathcal{U}$ implies the weak convergence of its components in $H^{1;a}(B^+_1)$. We know that $\int_{B^+_1}y^a|\nabla u|^2$ is weakly lower semicontinuous in $H^{1;a}(B^+_1)$ since it is the sum of the norm of the Hilbert space, which is weakly lower semicontinuous and of the $L^2(y^a\mathrm{d}z)$-norm, which is weakly continuous by Sobolev compact embeddings. Then, $\int_{\partial^0 B^+_1}\beta u^2v^2$ is weakly lower semicontinuous by the Fatou lemma; in fact, up to a subsequence, by the trace theorem, the weak convergence implies that $u_n\longrightarrow u$ and $v_n\longrightarrow v$ in $L^2(\partial B^+_1;\mathrm{d}\mu)$ where $\mathrm{d}\mu=y^a\mathrm{d}S(z)$ over $\partial^+B^+_1$ and $\mathrm{d}\mu=\mathrm{d}x$ over $\partial^0B^+_1$, and hence that $\beta u^2_n(z)v^2_n(z)\longrightarrow \beta u^2(z)v^2(z)$ for almost every $z\in\partial^0B^+_1$ with respect to the $2$-dimensional Lebesgue measure. So, we get the weak lower semicontinuity of $I$ as the sum of weakly lower semicontinuous pieces.

To show that $I$ is coercive, we want that
\begin{equation}\label{coer}
I(u,v)\geq\frac{1}{2}\int_{B^+_1}y^a(|\nabla u|^2+|\nabla v|^2)\longrightarrow+\infty,\qquad\mathrm{as}\quad||(u,v)||\longrightarrow+\infty,
\end{equation}
where $||(u,v)||^2=\int_{B^+_1}y^a(|\nabla u|^2+|\nabla v|^2+u^2+v^2)$. Recalling that $(u,v)=(G+u_0,H+v_0)\in\mathcal{U}$ where $(u_0,v_0)\in H^{1;a}_{\partial^+B^+_1}(B^+_1)\times H^{1;a}_{\partial^+B^+_1}(B^+_1)$ and that Poincar\'e inequality holds for such functions, then \eqref{coer} is a simple computation.

Hence, we have a nontrivial minimizer $(u,v)$ of $I$ over $\mathcal{U}$. Obviously also $(|u|,|v|)$ is a minimizer and hence we can assume that such a minimizer is nonnegative. Let us define the group $\mathbb{G}$ of all the reflections $\sigma_i$, with $i=1,...,k$ endowed with the composition between reflections. Let us define the action
\begin{eqnarray}\label{group2}
&& X\times \mathbb{G}\longrightarrow X\nonumber\\
&&[(u,v),g]\longmapsto (v\circ g,u\circ g).
\end{eqnarray}
That is, for $g=\sigma_i$, it holds $$[(u,v),\sigma_i]=(v\circ \sigma_i,u\circ \sigma_i),$$ and for $g=\sigma_i\circ\sigma_j$, it holds $$[(u,v),\sigma_i\circ\sigma_j]=[[(u,v),\sigma_i],\sigma_j]=[(v\circ \sigma_i,u\circ \sigma_i),\sigma_j]=(u\circ \sigma_i\circ\sigma_j,v\circ \sigma_i\circ\sigma_j).$$
It is easy to check that this action is isometric and that the functional $I$ is invariant with respect to this action. Since $\mathcal{U}=Fix(\mathbb{G})$, by the principle of symmetric criticality of Palais, the minimizer $(u,v)$ is also a nonnegative critical point for $I$ over the whole $X$ and hence a weak solution to \eqref{systb1} with the desired property 1.

Finally, using the fact that $(u,v)$ is a minimizer of $I$ in $\mathcal{U}$ and also that $(G,H)\in\mathcal{U}$, we get the condition 2; that is,
\begin{equation}\label{limitNabove}
I(u,v)\leq I(G,H)=\frac{1}{2}\int_{B^+_1}y^a(|\nabla G|^2+|\nabla H|^2)=\frac{d}{2}
\end{equation}
since $G$ and $H$ are segregated in $\partial^0B^+_1$ and are homogeneous of degree $d$. In \eqref{limitNabove} we have used \eqref{h=1} and the Euler formula for homogeneous functions.
\endproof

\subsection{Blow-up and uniform bounds on compact sets}
Let us consider the sequence of solutions $(u_\beta,v_\beta)$ constructed in Lemma \ref{teo3}. Thanks to the uniform bound given by condition 2, and the fact that the functional $I$ is coercive, we obtain uniform boundedness in $H^{1;a}(B^+_1)$ for both components of such a sequence. Hence, letting $\beta\longrightarrow+\infty$, there exists a weak limit $(U,V)$.

We remark that solutions $(u_\beta,v_\beta)$ of \eqref{systb1} are strictly positive in the open $B^+_1$ by maximum principles for $L_a$-subharmonic functions (see \cite{cabresire}), and for the same reason they are strictly positive also in $\partial^+B^+_1$ since it holds a maximum principle for $(g,h)$ over $S^2_+$. Moreover, they are strictly positive also in $\partial^0B^+_1$. By contradiction $u_{\beta}(z_0)=0$ for a point $z_0\in\partial^0B^+_1$ that is a minimum for $u_{\beta}$. By the Hopf lemma $\partial_y^au_{\beta}(z_0)<0$ (Proposition 4.11 in \cite{cabresire}) but the boundary condition imposed over the flat part of the boundary says that $-\partial_y^au_{\beta}(z_0)=u_{\beta}(z_0)v^2_{\beta}(z_0)=0$. Hence, they are able to assume value zero only on $S^1=\partial S^2_+$.

Moreover $(u_\beta,v_\beta)$ must attain their supremum in $\partial^+B^+_1$. Let us consider for example the component $u_{\beta}$. Its supremum must be attained by a point $z_0\in\partial B^+_1$ for the maximum principle but this point can not be on $\partial^0B^+_1$ by the Hopf lemma. In fact, we would obtain that $\partial_y^au_{\beta}(z_0)>0$ but
\begin{equation}
\partial_y^au_{\beta}=-\beta u_{\beta}v_{\beta}\leq 0\qquad\mathrm{in}\quad\partial^0 B^+_1,
\end{equation}
by boundary conditions and since $(u_{\beta},v_{\beta})$ are nonnegative.

So, all the functions $u_\beta$ are nonnegative, $L_a$-harmonic and such that
\begin{equation}\label{supppppp}
\sup_{\overline{B^+_1}}u_\beta=\sup_{\partial^+B^+_1}u_\beta=\sup_{\partial^+B^+_1}g=:A<+\infty.
\end{equation}
Moreover, thanks to \eqref{h=1}, $A>0$ since
\begin{equation}\label{lowerbound1}
1= \int_{\partial^+B^+_1}y^a(g^2+h^2)\leq 2\mu(\partial^+B^+_1)\left(\sup_{\partial^+B^+_1}g\right)^2=cA^2.
\end{equation}
The same holds for the functions $v_\beta$. Now, by this uniform boundedness obtained in $L^\infty(B^+_1)$, we can apply Theorem 1.1 in \cite{terverzil2}, obtaining for our sequence of solutions uniform boundedness in $C^{0,\alpha}_{loc}(\overline{B^+_1})$. This implies that the convergence of $(u_\beta,v_\beta)$ to $(U,V)$ is also uniform on every compact set in $B^+_1$. Moreover, since $A>0$, we get that the limit functions $(U,V)$ are not trivial and also nonnegative.

Likewise Soave and Zilio have done in \cite{soazil1} for the local case, we use a blow-up argument. For a radius $r_\beta\in(0,1)$ to be determined, we define
\begin{equation}\label{blowupbeta}
(\overline{u}_\beta,\overline{v}_\beta)(z):=\beta^{1/2}r_\beta^{s}(u_\beta,v_\beta)(r_\beta z).
\end{equation}
It is easy to check that such a blow-up sequence satisfies for every fixed $\beta>1$ the problem
\begin{equation}\label{systbrbeta}
\begin{cases}
L_a u=L_a v=0 & \mathrm{in}\quad B^+_{1/r_\beta}, \\
-\partial^a_y u=uv^2, \ -\partial^a_y v=vu^2 & \mathrm{in}\quad \partial^0B^+_{1/r_\beta}.
\end{cases}
\end{equation}
As in \cite{soazil1}, the choice of $r_\beta\in(0,1)$ is suggested by the following result.
\begin{Lemma}\label{lem4}
For any fixed $\beta>1$ there exists a unique $r_\beta\in(0,1)$ such that
\begin{equation}\label{Hoverubeta}
\int_{\partial^+B^+_1}y^a(\overline{u}_\beta^2+\overline{v}_\beta^2)=1.
\end{equation}
Moreover $r_\beta\longrightarrow 0$ as $\beta\longrightarrow+\infty$.
\end{Lemma} 
\proof
In order to prove \eqref{Hoverubeta}, we have to find for any fixed $\beta>1$, a radius $r_\beta\in(0,1)$ such that $\beta r^{2s}_\beta H((u_\beta,v_\beta),r_\beta)=1$. The strict increasing monotonicity of $r\mapsto H(r)$ (see e.g. \cite{terverzil,wangwei}) implies that also the function $r\mapsto\beta r^{2s}H((u_\beta,v_\beta),r)$ is strictly increasing and regular. Hence, for $\beta>1$ fixed,
\begin{equation}\label{conti1}
\lim_{r\to 0}\beta r^{2s}H((u_\beta,v_\beta),r)=\lim_{r\to 0}\beta r^{2s-2-a}\int_{\partial^+B^+_r}y^a(u_\beta^2+v_\beta^2)=\beta(u_\beta^2(0)+v_\beta^2(0))\lim_{r\to 0}r^{2s}=0.
\end{equation}
Moreover, by \eqref{h=1}, $\beta H((u_\beta,v_\beta),1)=\beta>1$. Obviously, existence and uniqueness of $r_\beta$ follow. If, seeking a contradiction, it would exist $\overline{r}>0$ such that for any $\beta>1$ it holds $r_\beta\geq\overline{r}$, then by the monotonicity recalled above and using \eqref{rappHH} and \eqref{h=1}, we get
\begin{equation}\label{conti2}
1=\beta r^{2s}_\beta H((u_\beta,v_\beta),r_\beta)\geq\beta \overline{r}^{2s} H((u_\beta,v_\beta),\overline{r})\geq c\beta\overline{r}^{2d+2s}H((u_\beta,v_\beta),1)=c\beta.
\end{equation}
So, we get a contradiction for choices of $\beta>1/c$.
\endproof

\subsubsection{Proof of Theorem \ref{tteo2}}
In order to prove Theorem \ref{tteo2}, we want to prove the existence of positive functions $(\overline{U},\overline{V})$ which solve \eqref{cafsilvess} and such that $(\overline{u}_\beta,\overline{v}_\beta)\longrightarrow(\overline{U},\overline{V})$ uniformly on compact sets of $\mathbb{R}^{3}_+$ with $N((\overline{U},\overline{V}),r)\leq d$ for any $r>0$. Hence, according to \cite{wangwei}, we would obtain in the case $N=2$ a solution of \eqref{cafsilvess} which grows asymptotically no more than
\begin{equation}\label{dkgrowth}
\overline{U}(x,y)+\overline{V}(x,y)\leq c\left(1+|x|^2+y^2\right)^{d/2},
\end{equation}
with $d=d(k)\in[s,2s)$. Moreover, we will prove that the growth rate of this solution is exactly equal to $d$.

Thanks to the monotonicity of the frequency and conditions \eqref{h=1} and \eqref{limitNabove}, we get for any $\beta>1$ and $r\in(0,1/r_\beta)$,
\begin{equation}\label{conti3}
N((\overline{u}_\beta,\overline{v}_\beta),r)\leq N((\overline{u}_\beta,\overline{v}_\beta),1/r_\beta)=\frac{2I(u_\beta,v_\beta)}{H((u_\beta,v_\beta),1)}\leq d.
\end{equation}
Moreover, for any $\beta>1$ large, for any $1\leq r\leq\frac{1}{r_\beta}$, using \eqref{rappHH}, we obtain the following upper bound which does not depend on $\beta$,
\begin{equation}\label{conti4}
H((\overline{u}_\beta,\overline{v}_\beta),r)\leq  H((\overline{u}_\beta,\overline{v}_\beta),1)e^{\frac{d}{1-a}}r^{2d}=e^{\frac{d}{1-a}}r^{2d}.
\end{equation}
Since for every $\beta>0$ the functions $(\overline{u}_\beta,\overline{v}_\beta)$ have $-\partial^a_y\overline{u}_\beta\geq 0$, $-\partial^a_y\overline{v}_\beta\geq 0$, then their extensions to $B_{1/r_\beta}$ (through even reflections with respect to $\{y=0\}$) satisfy the requirements of Lemma A.2 in \cite{wangwei}. Then it holds that both the components $\overline{u}_\beta$ and $\overline{v}_\beta$ satisfy
\begin{equation}\label{lemA2}
\sup_{B^+_r}u\leq c\left(\frac{1}{r^{3+a}}\int_{B^+_{2r}}y^au^2\right)^{1/2}.
\end{equation}
Hence, using \eqref{conti3}, \eqref{conti4} and \eqref{lemA2}, we get the upper bound
\begin{equation}\label{Himplicasup}
\left(\sup_{B^+_r}(\overline{u}_\beta+\overline{v}_\beta)\right)^2\leq C(r)H((\overline{u}_\beta,\overline{v}_\beta),r)\leq C(r);
\end{equation}
that is both components of the sequence $(\overline{u}_\beta,\overline{v}_\beta)$ are uniformly bounded in $L^\infty(B^+_r)$, independently from $\beta$ large enough. This gives us uniform boundedness in $C^{0,\alpha}_{loc}(\overline{B^+_r})$ (see \cite{terverzil}) and so, up to consider a subsequence, this ensures the convergence to a nontrivial nonnegative function on compact subsets of $\overline{B^+_r}$. By the arbitrariness of the choice of $r\geq 1$ done, we obtain such a convergence on every compact set in $\overline{\mathbb{R}^{3}_+}$. Since for $\beta\longrightarrow+\infty$ we have $1/r_{\beta}\longrightarrow+\infty$, then the limit $(\overline{U},\overline{V})$ is a nonnegative solution to \eqref{cafsilvess} with $N((\overline{U},\overline{V}),r)\leq d$ for any $r>0$ using the uniform convergence and \eqref{conti3}. Hence \eqref{dkgrowth} follows.

Now, we have to verify that $(\overline{U},\overline{V})$ are strictly positive in $\overline{\mathbb{R}^{3}_+}$. Obviously, by construction they are nonnegative in $\overline{\mathbb{R}^{3}_+}$ and strictly positive in $\mathbb{R}^{3}_+$ by maximum principles. Moreover, it is impossible that one component has a zero in $\partial\mathbb{R}^{3}_+$. By contradiction let $z_0\in\partial\mathbb{R}^{3}_+$ be such that $\overline{U}(z_0)=0$. By the Hopf lemma it would be $\partial_y^a\overline{U}(z_0)<0$ since this point is a minimum. But, by the boundary condition we get the contradiction
\begin{equation}
-\partial_y^a\overline{U}(z_0)=\overline{U}(z_0)\overline{V}^2(z_0)=0.
\end{equation}
Hence, we want to show that the asymptotic growth rate is exactly equal to $d=d(k)$. Seeking a contradiction, let $N((\overline{U},\overline{V}),r)\leq d(k)-\varepsilon$ for any $r>0$. By the Almgren monotonicity formula, there exists the limit $\lim_{r\to+\infty}N((\overline{U},\overline{V}),r):=\overline{d}\leq d(k)-\varepsilon$. We replicate the blow-down construction performed in section 2.3 on the solution $(\overline{U},\overline{V})$, obtaining the convergence in $C^{0,\alpha}_{loc}(\overline{\mathbb{R}^{3}_+})$ of the blow-down sequence to a couple of $\overline{d}$-homogeneous functions segregated in $\partial\mathbb{R}^{3}_+$. The spherical parts of this functions are eigenfunctions with same eigenvalue $\overline{\lambda}$ of the Laplace-Beltrami-type operator on $S^2_+$ which inherit their symmetries from the functions $(u_\beta,v_\beta)$ (see \eqref{symmmm}). In fact such symmetries hold also for the blow-up sequence $(\overline{u}_\beta,\overline{v}_\beta)$ constructed in \eqref{blowupbeta} and hence also for $(\overline{U},\overline{V})$, thanks to the uniform convergence on compact sets. By the condition $\overline{d}\leq d(k)-\varepsilon$ over the characteristic exponent, hence we have $\overline{\lambda}<\lambda_1^s(k)$ using \eqref{lambdad}, but by the minimality of $\lambda_1^s(k)$ we would have $\overline{\lambda}\geq\lambda_1^s(k)$ since its eigenfunction is a competitor for the problem defined in \eqref{InfInfk}, and hence we get the contradiction.

Eventually, let us say that these prescribed growth solutions for \eqref{cafsilvess} in space dimension $N=2$ are also solutions with the same properties for the same problem in any higher dimension.
This concludes the proof of Theorem \ref{tteo2}.

\section{Multidimensional entire solutions}
In this section we will show the existence of $N$-dimensional entire solutions to \eqref{cafsilvess} which can not be obtained by adding coordinates in a constant way starting from a $2$-dimensional solution. Actually, we will establish a more general result for system \eqref{cafsilvess} in case of $k$-component; that is, considering solutions $\mathbf{u}:=(u_1,...,u_k)$ to
\begin{equation}\label{cafsilvesss}
\begin{cases}
L_au_i=0, &\mathrm{in} \ \mathbb{R}^{N+1}_+,\\
u_i>0, &\mathrm{in} \ \overline{\mathbb{R}^{N+1}_+},\\
-\partial^a_yu_i=u_i\sum_{j\neq i}u_j^2, &\mathrm{in} \ \partial\mathbb{R}^{N+1}_+,
\end{cases}
\end{equation}
for any $i=1,...,k$. In what follows, we adapt the results for the local case in \cite{soazil1} to the fractional setting.

First of all, we remark that also in the case of $k$-components Theorem \ref{tteo1} holds; that is, solutions to \eqref{cafsilvesss} have a universal bound on the growth rate at infinity given by
\begin{equation}\label{ks}
u_1(x,y)+...+u_k(x,y)\leq c(1+|x|^2+y^2)^{s}.
\end{equation}
In fact, also in this setting a Pohozaev inequality holds (see \cite{terverzil2}); that is, for any $x_0\in\mathbb{R}^N$ and $r>0$,
\begin{eqnarray}
(N-1+a)\int_{B_r^+(x_0,0)}y^a\sum_{i=1}^k|\nabla u_i|^2&=&r\int_{\partial^+B_r^+(x_0,0)}y^a\sum_{i=1}^k|\nabla u_i|^2-2y^a\sum_{i=1}^k|\partial_r u_i|^2\nonumber\\
&+&r\int_{S^{N-1}_r(x_0,0)}\sum_{i,j<i}u_i^2u_j^2-N\int_{\partial^0B_r^+(x_0,0)}\sum_{i,j<i}u_i^2u_j^2.
\end{eqnarray}
Moreover, let us recall the following definitions
\begin{equation}\label{E}
E(r,x_0;\mathbf{u}):=\frac{1}{r^{N-1+a}}\left(\int_{B_r^+(x_0,0)}y^a\sum_{i=1}^k|\nabla u_i|^2 \ + \ \int_{\partial^0B_r^+(x_0,0)}\sum_{i,j<i}u_i^2u_j^2\right),
\end{equation}
and
\begin{equation}\label{H}
H(r,x_0;\mathbf{u}):=\frac{1}{r^{N+a}}\int_{\partial^+B_r^+(x_0,0)}y^a\sum_{i=1}^ku_i^2.
\end{equation}
Hence, defining the frequency as $N(r,x_0;\mathbf{u}):=\frac{E(r,x_0;\mathbf{u})}{H(r,x_0;\mathbf{u})}$, the Almgren monotonicity formula holds; that is, the frequency $N(r,x_0;\mathbf{u})$ is non decreasing in $r>0$ (the proof is as in \cite{wangwei}). Since the bound \eqref{dgrowth} found in \cite{wangwei} also holds in the case of solutions to \eqref{cafsilvesss}, one can apply the procedure seen in the proof of Theorem \ref{tteo1} obtaining eventually \eqref{ks}.

Let us denote by $\mathcal{O}(N)$ the orthogonal group of $\mathbb{R}^N$ and by $\mathfrak{G}_k$ the symmetric group of permutations of $\{1,...,k\}$. We assume the existence of a homomorphism $h:\mathbb{G}<\mathcal{O}(N)\to\mathfrak{G}_k$ with $\mathbb{G}$ a nontrivial subgroup. Hence, let us define the equivariant action of $\mathbb{G}$ on $H^{1;a}(\mathbb{R}^{N+1}_+,\mathbb{R}^k)$ so that
\begin{eqnarray}\label{groupaction}
&& H^{1;a}(\mathbb{R}^{N+1}_+,\mathbb{R}^k)\times \mathbb{G}\longrightarrow H^{1;a}(\mathbb{R}^{N+1}_+,\mathbb{R}^k)\nonumber\\
&&[\mathbf{u},g]\longmapsto (u_{(h(g))^{-1}(1)}\circ g,...,u_{(h(g))^{-1}(k)}\circ g),
\end{eqnarray}
where $\circ$ denotes the usual composition of functions. Let us define the space of the $(\mathbb{G},h)$-equivariant functions as
\begin{equation}
H_{(\mathbb{G},h)}:=Fix(\mathbb{G})=\{\mathbf{u}\in H^{1;a}(\mathbb{R}^{N+1}_+,\mathbb{R}^k) \ : \ \mathbf{u}\circ g=\mathbf{u} \ \forall g\in \mathbb{G}\}.
\end{equation}
As in \cite{soazil1}, we give the following definition.
\begin{Definition}\label{defq}
Let $k\in\mathbb{N}$, $\mathbb{G}<\mathcal{O}(N)$ be a nontrivial subgroup and $h:\mathbb{G}\to\mathfrak{G}$ a homomorphism. We say that the triplet $(k,\mathbb{G},h)$ is admissible if there exists $\mathbf{u}\in H_{(\mathbb{G},h)}$ such that
\begin{enumerate}
\item[$(i)$] $u_i\geq 0$ and $u_i\neq 0$ in $\mathbb{R}^{N+1}_+$ for any $i=1,...,k$,
\item[$(ii)$] $u_iu_j=0$ in $\mathbb{R}^N$ for any $i,j=1,...,k$ with $i\neq j$,
\item[$(iii)$] there exist $g_2,...,g_k\in\mathbb{G}$ such that $u_i=u_1\circ g_i$ for any $i=2,...,k$.
\end{enumerate}
\end{Definition}
We remark that if the triplet $(k,\mathbb{G},h)$ is admissible, then all the $(\mathbb{G},h)$-equivariant functions satisfy $(iii)$ of the definition with the same elements $g_2,...,g_k$. Moreover it holds $(h(g_i))^{-1}(i)=1$ for any $i=1,...,k$, and hence equivariant functions satisfy
\begin{equation}
u_i=u_{(h(g_i))^{-1}(i)}\circ g_i=u_1\circ g_i;
\end{equation}
that is, if the triplet is admissible, then any equivariant function $\mathbf{u}$ is determined by its first component $u_1$ and by knowing the elements $g_2,...,g_k$.

\subsection{Optimal $k$-partition problem}
Let us define the set of $k$-partitions of $S^{N-1}$ as
\begin{eqnarray}\label{partk}
\mathcal{P}^k &=&\{(\omega_1,...,\omega_k) : \ \omega_i\subset S^{N-1} \mathrm{ \ open}, \ \omega_i\cap\omega_j=\emptyset, \ \bigcup_{i=1}^k\overline{\omega_i}= S^{N-1}, \ \overline{\omega_i}\cap\overline{\omega_j} \ \mathrm{is \ a}\nonumber\\
&&(N-2)-\mathrm{dimensional \ smooth \ submanifold, \ }\forall i,j=1,...,k, \ j\neq i\}.
\end{eqnarray}
Let $(k,\mathbb{G},h)$ be an admissible triplet. We denote by
\begin{eqnarray}
\Lambda_{(\mathbb{G},h)}&=&\{\mathbf{u}\in H^{1;a}(S^N_+,\mathbb{R}^k) : \ \mathbf{u} \ \mathrm{is \ the \ restriction \ to \ } S^N_+ \ \mathrm{of \ a \ } (\mathbb{G},h)-\mathrm{equivariant \ function}\nonumber\\
&& \mathrm{fulfilling \ } (i),(ii),(iii) \ \mathrm{in \ Definition \ } \ref{defq}, \ \mathrm{such \ that \ } \int_{S^N_+}y^au_i^2=1 \ \forall i=1,...,k \}.
\end{eqnarray}
Obviously, assuming that the triplet is admissible, up to consider a normalization of the components in $L^2(S^N_+;\mathrm{d}\mu)$, it follows that $\Lambda_{(\mathbb{G},h)}$ is not empty. Moreover, by conditions $(i)$ and $(ii)$ one has that for any element $\mathbf{u}\in\Lambda_{(\mathbb{G},h)}$ there exists a $k$-partition $(\omega_1,...,\omega_k)\in\mathcal{P}^k$ such  that $u_i=0$ in $S^{N-1}\setminus\omega_i$ for any $i=1,...,k$. Let us consider the following minimization problem
\begin{equation}\label{mink}
\inf_{\mathbf{u}\in\Lambda_{(\mathbb{G},h)}}I(\mathbf{u}),
\end{equation}
where
\begin{equation}
I(\mathbf{u})=\frac{1}{2}\int_{S^N_+}y^a\sum_{i=1}^k|\nabla_{S^N}u_i|^2.
\end{equation}
One can easily check that problem \eqref{mink} produces a nontrivial nonnegative minimizer $\mathbf{u}$ in $\Lambda_{(\mathbb{G},h)}$, and since the functional $I$ is invariant with respect to the action in \eqref{groupaction}, applying the principle of criticality of Palais, we obtain that such a minimizer is also a solution to an eigenvalue problem; that is, its components satisfy for any $i,j=1,...,k$, $j\neq i$
\begin{equation}
\begin{cases}
-L_a^{S^N}u_i=y^a\lambda u_i, &\mathrm{in} \ S^{N}_+,\\
u_iu_j=0, &\mathrm{in} \ S^{N-1},\\
u_i\partial^a_yu_i=0, &\mathrm{in} \ S^{N-1},
\end{cases}
\end{equation}
where $\lambda=\int_{S^N_+}y^a|\nabla_{S^N}u_1|^2=...=\int_{S^N_+}y^a|\nabla_{S^N}u_k|^2$, by condition $(iii)$ and the invariance of $I$ with respect to the group action. Moreover there exists a $k$-partition $(\omega_1,...,\omega_k)\in\mathcal{P}^k$ such that for any $i=1,...,k$ it holds $u_i=0$ in $S^{N-1}\setminus\omega_i$. We want to prove the $C^{0,\alpha}(\overline{S^N_+})$-regularity for the components of $\mathbf{u}$ via the convergence of solutions of $\beta$-problems over $S^N_+$ to our eigenfunctions.
Let us now consider the following set of functions
\begin{eqnarray}\label{Ckeq}
\Gamma_{(\mathbb{G},h)}&=&\{\mathbf{u}\in H^{1;a}(S^N_+,\mathbb{R}^k) : \ \mathbf{u} \ \mathrm{is \ the \ restriction \ to \ } S^N_+ \ \mathrm{of \ a \ } (\mathbb{G},h)-\mathrm{equivariant \ function}\nonumber\\
&& \mathrm{fulfilling \ } (i),(iii) \ \mathrm{in \ Definition \ } \ref{defq}, \ \mathrm{such \ that \ } \int_{S^N_+}y^au_i^2=1 \ \forall i=1,...,k \}.
\end{eqnarray}
This space is trivially not empty since $\Lambda_{(\mathbb{G},h)}\subseteq\Gamma_{(\mathbb{G},h)}$.

Hence, for any $\beta>0$, we consider the following minimization problem
\begin{equation}\label{minbetalambdaeq}
\inf_{\mathbf{u}\in\Gamma_{(\mathbb{G},h)}}J_\beta(\mathbf{u}),
\end{equation}
with
\begin{equation}
J_\beta(\mathbf{u})=\frac{1}{2}\int_{S^N_+}y^a\sum_{i=1}^k|\nabla_{S^N}u_i|^2 \ + \ \frac{1}{2}\int_{S^{N-1}}\beta \sum_{i<j}u_i^2u_j^2=I(\mathbf{u}) \ + \ \frac{1}{2}\int_{S^{N-1}}\beta \sum_{i<j}u_i^2u_j^2.
\end{equation}
It is easy to check that, for every $\beta>0$ fix, there exists a nonnegative minimizer $\mathbf{u}_\beta\in\Gamma_{(\mathbb{G},h)}$. Moreover, since $J_\beta$ is invariant with respect to the action in \eqref{groupaction}, we get that this minimizer is also a weak solution to the system
\begin{equation}\label{asd2eq}
\begin{cases}
-L_a^{S^N} u_{\beta,i}=y^a\lambda_{\beta,i} u_{\beta,i}, & \mathrm{in}\quad S^N_+, \\
-\partial^a_y u_{\beta,i}=\beta u_{\beta,i}\sum_{j\neq i}u^2_{\beta,j}, & \mathrm{in}\quad S^{N-1},
\end{cases}
\end{equation}
for any $i=1,...,k$, where $\lambda_{\beta,i}=\int_{S^N_+}y^a|\nabla_{S^N}u_{\beta,i}|^2+\int_{S^{N-1}}\beta u_{\beta,i}^2\sum_{j\neq i}u_{\beta,j}^2$.
Moreover, since the minimizer $\mathbf{u}\in\Lambda_{(\mathbb{G},h)}\subseteq\Gamma_{(\mathbb{G},h)}$, it holds that for any $\beta>0$, we get the uniform bound
\begin{equation}\label{boueq}
0\leq\frac{1}{2}\sum_{i=1}^k\lambda_{\beta,i}\leq 2J_\beta(\mathbf{u}_\beta)\leq 2J_\beta(\mathbf{u})=k\lambda.
\end{equation}
This uniform bound gives the weak convergence in $H^{1;a}(S^N_+;\mathbb{R}^k)$ of the $\beta$-sequence to a function $\mathbf{u}_\infty$ (any component has the same norm $\int_{S^N}y^a|\nabla_{S^N}u_{\infty,i}|=\lambda_\infty$). Moreover, since solutions to \eqref{asd2eq} are bounded in $C^{0,\alpha}(\overline{S^N_+})$ uniformly in $\beta>0$ for $\alpha>0$ small, as it is proved in \cite{terverzil2}, we obtain, up to consider a subsequence as $\beta\to+\infty$, that the convergence is uniform on compact sets and so that the limit satisfies the symmetry conditions. Moreover it holds that $\int_{S^{N-1}}\beta u_{\beta,i}^2u_{\beta,j}^2\to 0$ for any $i,j=1,...,k$ with $j\neq i$ (see Lemma 4.6 in \cite{terverzil2} and Lemma 5.6 in \cite{terverzil} for the details in the case $s=1/2$). So, the limit should have the components segregated on $S^{N-1}$; that is, $\mathbf{u}_\infty\in\Lambda_{(\mathbb{G},h)}$. Moreover, we have
\begin{equation}\label{bou2eq}
0\leq\frac{1}{2}\sum_{i=1}^k\lambda_{\beta,i}= J_\beta(\mathbf{u}_\beta)+\frac{1}{2}\int_{S^{N-1}}\beta \sum_{i<j}u_i^2u_j^2\leq \frac{k}{2}\lambda+\frac{1}{2}\int_{S^{N-1}}\beta \sum_{i<j}u_i^2u_j^2,
\end{equation}
and since for any $i=1,...,k$ one has $\lambda_{\beta,i}\to\lambda_\infty$, by \eqref{bou2eq} it follows that $\lambda_\infty\leq\lambda$. But by the minimality of $\mathbf{u}$ in $\Lambda_{(\mathbb{G},h)}$ we have also $\lambda\leq\lambda_\infty$, and hence we obtain that $\mathbf{u}_\infty$ and $\mathbf{u}$ own the same norm in $H^{1;a}(S^N_+;\mathbb{R}^k)$, and hence we can choose as a minimizer $\mathbf{u}_\infty$ which inherits the H\"older regularity up to the boundary.

\subsection{$(\mathbb{G},h)$-equivariant solutions}
In order to construct $(\mathbb{G},h)$-equivariant entire solutions to \eqref{cafsilvesss}, one can follow the construction given in section 3.7 and 3.8. Let us summarize the main steps: first, we construct $(\mathbb{G},h)$-equivariant $\beta$-solutions $\mathbf{u}_\beta$ on the unit half ball $B^+_1$; that is, solutions inheriting the symmetries given by the triplet $(k,\mathbb{G},h)$ and so that the boundary value on $\partial^+B_1^+$ is the minimizer $\mathbf{u}$ previously found (the proof follows from Lemma \ref{teo3}). Since any component $u_i$ of $\mathbf{u}$ has the same energy $\int_{S^N_+}y^a|\nabla_{S^N}u_i|^2=\lambda$, we can define the $d$-homogeneous extension of $\mathbf{u}$ to $\mathbb{R}^{N+1}_+$, where $d=d(\lambda)$; that is, $\overline{\mathbf{u}}=|z|^d\mathbf{u}(\frac{z}{|z|})$. This function gives a bound over the energy of our $\beta$-solutions; that is,
\begin{equation}
2F_\beta(\mathbf{u}_\beta)=\int_{B^+_1}y^a\sum_{i=1}^k|\nabla u_{i,\beta}|^2 \ + \ \beta\int_{\partial^0B^+_1}\sum_{i<j}u_{i,\beta}^2u_{j,\beta}^2\leq d.
\end{equation}
Hence, after rescaling (the right choice is given by an analogous of Lemma \ref{lem4}), by the blow-up argument, we get convergence to a positive $(\mathbb{G},h)$-equivariant entire solution $\mathbf{U}$ to \eqref{cafsilvesss} as $\beta\to+\infty$ on compact subsets of $\mathbb{R}^{N+1}_+$. Moreover, for any $r>0$, we get a bound over the Almgren frequency given by
\begin{equation}\label{almbf}
N(r;\mathbf{U})\leq d.
\end{equation}

\subsection{An admissible triplet $(2,\mathbb{G},h)$}
To conclude this section, we want to provide the existence, for simplicity in the case of two components, of multidimensional entire solutions to \eqref{cafsilvess} in $\mathbb{R}^{N+1}_+$ with $N\geq 3$ and such that they can not be obtained by adding coordinates in a constant way starting from a $2$-dimensional solution. Let $k=2$ and $\mathbb{G}<\mathcal{O}(N)$ be the nontrivial subgroup of symmetries generated by the reflections $\sigma_i$ with respect to the hyperplanes $\Sigma_i=\{x_i=0\}$ for any $i=1,...,N$. Let also $h:\mathbb{G}\to\mathfrak{G}_2$ be defined on the generators of $\mathbb{G}$ by $h(\sigma_i)=(1 \ \ 2)$ for every $i=1,...,N$ (the expression $(1 \ \ 2)$ denotes the cycle mapping 1 in 2 and 2 in 1). Let us consider the fundamental domain defined as the set $\mathcal{D}(2,\mathbb{G},h)=S^N_+\cap\{z=(x,y)\in\mathbb{R}^{N+1}_+: \ x_2>0,x_3>0,...,x_N>0\}$. Obviously there exists a couple of nontrivial and nonnegative functions $(f_1,f_2)$ such that $f_1\in H^{1;a}_0(\mathcal{D}(2,\mathbb{G},h)\cap\{x_1>0\})$ and  $f_2\in H^{1;a}_0(\mathcal{D}(2,\mathbb{G},h)\cap\{x_1<0\})$. Let us merge them in a unique function $v_1$ over the fundamental domain, and extend it to the whole of the hemisphere $S^N_+$ following the condition $v_1(z)=v_1(\sigma_i(\sigma_j(z)))$ for any $i,j=1,...,N$ (the values of $u_1$ over the fundamental domain are enough to define it on the hemisphere). In the same way, we can define the function $v_2$ so that $v_2(z)=v_1(\sigma_i(z))$ for every $i=1,...,N$. Let us normalize the two functions in $L^2(S^N_+;\mathrm{d}\mu)$. Let us also define the number $\nu=\int_{S^N_+}y^a|\nabla_{S^N}v_1|^2=\int_{S^N_+}y^a|\nabla_{S^N}v_2|^2$. The $d(\nu)$-homogeneous extension of $\mathbf{v}=(v_1,v_2)$ to $\mathbb{R}^{N+1}_+$ (the characteristic exponent is defined in \eqref{mono}) is an element of $H_{(\mathbb{G},h)}$ satisfying conditions $(i)$, $(ii)$ and $(iii)$, and hence, as a consequence, the triplet $(2,\mathbb{G},h)$ turns out to be admissible.

Hence, it is possible to apply the construction seen in the first part of this section, in order to construct a $(\mathbb{G},h)$-equivariant solution to \eqref{cafsilvess} depending on the minimizer of the problem \eqref{mink}. We want to show that it holds
\begin{equation}\label{minkq}
\inf_{\mathbf{u}\in\Lambda_{(\mathbb{G},h)}}I(\mathbf{u})=\lambda<\lambda_1^s(\emptyset).
\end{equation}
Let us define the set of 2-partitions
\begin{eqnarray}\label{part2}
\mathcal{P}^2_N &=&\{(\omega_1,\omega_2) : \ \omega_i\subset S^{N-1} \mathrm{ \ open}, \ \omega_1\cap\omega_2=\emptyset, \ \overline{\omega_1}\cup\overline{\omega_2}= S^{N-1},  \ \overline{\omega_1}\cap\overline{\omega_2} \ \mathrm{is \ a}\nonumber\\
&&(N-2)-\mathrm{dimensional \ smooth \ submanifold, \ }\ \mathcal{H}^{N-1}(\omega_1)=\mathcal{H}^{N-1}(\omega_2),\nonumber\\
&& z\in\omega_1\Leftrightarrow \sigma_i(z)\in\omega_2 \ \forall i=1,...,N\}.
\end{eqnarray}
Let us also introduce, for any element $(\omega_1,\omega_2)\in\mathcal{P}^2_N$ the space
\begin{equation}
\Lambda_{(\mathbb{G},h)}(\omega_1,\omega_2)=\{\mathbf{u}\in\Lambda_{(\mathbb{G},h)} : \ u_i=0 \ \mathrm{in \ } S^{N-1}\setminus\omega_i, \ \forall i=1,2\}.
\end{equation}
We remark that for any 2-partition, this space is not empty since the function $\mathbf{v}$ previously constructed is
contained. The minimization problem in \eqref{mink} can be expressed as
\begin{equation}
\inf_{\mathbf{u}\in\Lambda_{(\mathbb{G},h)}}I(\mathbf{u})=\inf_{(\omega_1,\omega_2)\in\mathcal{P}^2_N}\inf_{\mathbf{u}\in\Lambda_{(\mathbb{G},h)}(\omega_1,\omega_2)}I(\mathbf{u}).
\end{equation}
Let us consider the particular 2-partition $(\omega_1^N,\omega_2^N)\in\mathcal{P}^2_N$ so that $\omega_1^N\supset S^{N-1}\cap\{x_1>0,...,x_N>0\}$. Obviously one has
\begin{equation}
\lambda\leq\inf_{\mathbf{u}\in\Lambda_{(\mathbb{G},h)}(\omega_1^N,\omega_2^N)}I(\mathbf{u})=\lambda_N.
\end{equation}
Therefore, by considerations over the symmetries and \eqref{coloper}, it is easy to see that
\begin{equation}
\lambda_N=\frac{\lambda_1^s(\omega_1^N)+\lambda_1^s(\omega_2^N)}{2}=\lambda_1^s(\omega_1^N)<\lambda_1^s(\emptyset),
\end{equation}
since $\mathcal{H}^{N-1}(\omega_1^N)>0$. Hence, \eqref{minkq} holds.

Now we want to show that the equivariant entire solution $\mathbf{U}=(U_1,U_2)$ obtained depends on any $x_i$-variable for any $i=1,...,N$. Thanks to the bound in \eqref{almbf} and the condition \eqref{minkq}, we get that
\begin{equation}\label{almbf2}
N(r;\mathbf{U})\leq d<d(\lambda_1^s(\emptyset))=2s.
\end{equation}
Let us suppose by contradiction that $\mathbf{U}$ does not depend on the variable $x_1$ (we can choose it without loss of generality). Then, considering the reflection $\sigma_1$, one has for any $z\in\overline{\mathbb{R}^{N+1}_+}$
\begin{equation}\label{=zz}
U_1(z)=U_1(\sigma_1(z))=U_2(z).
\end{equation}
Let us proceed now by a blow-down construction as in section 2.3 for the proof of Theorem \ref{tteo1}. The limit of the blow-down sequence is a couple $(u_\infty,v_\infty)$ of functions solving
\begin{equation}\label{sblowcseg2}
\begin{cases}
L_a u_{\infty}=L_a v_{\infty}=0 & \mathrm{in}\quad\mathbb{R}^{N+1}_+, \\
u_{\infty}\partial^a_yu_{\infty}=v_{\infty}\partial^a_yv_{\infty}=0 & \mathrm{in}\quad\partial\mathbb{R}^{N+1}_+, \\
u_{\infty}v_{\infty}=0 & \mathrm{in}\quad\partial\mathbb{R}^{N+1}_+.
\end{cases}
\end{equation}
By the uniform convergence, condition \eqref{=zz} says that $u_\infty=v_\infty$ in $\overline{\mathbb{R}^{N+1}_+}$, and by the segregation condition also that $u_\infty=v_\infty=0$ in $\partial\mathbb{R}^{N+1}_+$. Moreover, such solutions have the form
\begin{equation*}\label{uinf2}
u_{\infty}(r,\theta)=v_{\infty}(r,\theta)=r^dg(\theta),
\end{equation*}
where $g$ is defined on the upper hemisphere $S^N_+=\partial^+B^+_1$. Since we have constructed the blow-down sequence so that $H(1;\mathbf{U}_R)=1$, then
\begin{equation}\label{xxx}
\int_{S^N_+}y^ag^2=1/2.
\end{equation}
Since $d<2s$, we can apply a Liouville type result (see Proposition 3.1 in \cite{terverzil2}) in order to conclude that $u_\infty$ and $v_\infty$ should be trivial everywhere, in contradiction with condition \eqref{xxx}.

\end{document}